# A CHARACTERIZATION OF THE OPTIMAL RISK-SENSITIVE AVERAGE COST IN FINITE CONTROLLED MARKOV CHAINS[1]

By Rolando Cavazos-Cadena and Daniel Hernández-Hernández

*Universidad Autónoma Agraria Antonio Narro
and Centro de Investigación en Matemáticas*

This work concerns controlled Markov chains with finite state and action spaces. The transition law satisfies the simultaneous Doeblin condition, and the performance of a control policy is measured by the (long-run) risk-sensitive average cost criterion associated to a positive, but otherwise arbitrary, risk sensitivity coefficient. Within this context, the optimal risk-sensitive average cost is characterized via a minimization problem in a finite-dimensional Euclidean space.

**1. Introduction.** This work concerns discrete-time Markov decision processes (MDPs), where the controller selects actions from a finite set, and the corresponding controlled process takes values on a finite set $S$. The decision maker is supposed to be risk-averse with constant risk sensitivity coefficient $\lambda > 0$, and the performance index of a control policy is measured by the (long-run) risk-sensitive average cost criterion. Under the simultaneous Doeblin condition in Assumption 2.1, the main result of the paper, stated as Theorem 3.5, provides a characterization of the optimal value function $J^*(\lambda, \cdot)$ for arbitrary $\lambda > 0$. Roughly, this theorem shows that the optimal value function is the infimum of a family $\mathcal{G}$ of functions on the state space, a conclusion that, as described in the following section, is similar to results already available for classical risk-neutral criteria. However, at the same time this characterization reflects an interesting and important contrast with the risk-neutral average cost index which is illustrated in Example 2.2, namely, when $\lambda$ is large enough, the costs incurred while the system stays at transient states, which can be visited only at "early stages" of the decision process,

Received April 2002; revised December 2003.
[1]Supported in part by PSF Organization Grant 100-100/02–03 and Conacyt Grant 37643-E.

have a definite impact in the risk-sensitive average performance criterion. This feature implies that, even when the Markov chain associated with each stationary policy has a single recurrent class, the risk-sensitive optimal average cost is not necessarily constant, and that in this case the optimality equation may have no solution at all. Such a potentially complex behavior of $J^*(\lambda, \cdot)$ when $\lambda > 0$ is unrestricted is actually covered by the characterization in Theorem 3.5, and highlights the main difference between the results in this paper and those already available, which concern the case in which $\lambda$ is small enough to guarantee that $J^*(\lambda, \cdot)$ is constant and its value is determined via the optimality equation; see, for instance, [3] or [14] for the discrete case, or [8] for MDPs over Borel spaces.

The study of MDPs endowed with the risk-sensitive average criterion can be traced back, at least, to the seminal work of Howard and Matheson [17], where models with finite state and action spaces were studied assuming the following condition (C): Under each stationary policy the whole state space is an aperiodic communicating class. In this context, the Perron–Frobenius theory of positive matrices [7] was used to show that, for every $\lambda > 0$, the $\lambda$-sensitive average cost associated to each stationary policy is a constant function, and its value $\gamma$ can be characterized via the corresponding Poisson equation; see also [11]. The Perron–Frobenius theory provides also a link between risk-sensitive control and the Donsker–Varadhan theory of large deviations [9]. It is well known that, under suitable recurrence conditions, the occupation measure of a Markov process satisfies the large deviation principle, with rate function given by the convex conjugate of a long-run expected rate of exponential growth function. It is also worth mention that some optimal investment models can be formulated as risk-sensitive control problems, for assets dynamics models affected by economic factors, where the goal is to maximize the growth rate of the expected utility of wealth [1, 2, 12]. This kind of problems are also linked with the deterministic model of optimal economic development proposed by Gale and Neumann [10, 13].

The organization of the paper is as follows. In Section 2 a formal description of the model is presented, the potentially complex dependence of $J^*(\lambda, \cdot)$ on $\lambda > 0$ is explicitly shown and, after describing the main theorem, an outline of the strategy that will be used to prove the characterization result is given. In Section 3 a fundamental min–max equation satisfied by the optimal value function is established, and such an equality is used as one of the conditions in the definition of the family $\mathcal{G}$ in terms of which $J^*(\lambda, \cdot)$ is characterized in Theorem 3.5. After identifying the difficulties in proving this result, the necessary technical preliminaries are established in Sections 4–6 and, finally, the main theorem is proved in Section 7.

NOTATION.    Throughout the remainder $\mathbb{R}$ and $\mathbb{N}$ stand for the set of real numbers and nonnegative integers, respectively. Given a finite set $S$, the space of all real-valued functions defined on $S$ is denoted by $\mathcal{B}(S)$, and for each $C \in \mathcal{B}(S)$

$$\|C\| := \max_{w \in S} |C(w)|$$

is the corresponding maximum norm. The indicator function associated to an event $W$ is denoted by $I[W]$ and, even without explicit reference, all relations involving conditional expectations are supposed to hold almost surely with respect to the underlying probability measure.

**2. Decision model and outline of the work.** Let an MDP be specified by $M = \langle S, A, \{A(x)\}, C, P\rangle$ where the state space $S$ and the action set $A$ are finite sets endowed with the discrete topology and, for each $x \in S$, $A(x) \subset A$ is the nonempty subset of admissible actions at state $x$; the set $\mathbb{K}$ of admissible pairs is defined by $\mathbb{K} := \{(x,a) | a \in A(x), x \in S\}$ and is considered as a topological subspace of $S \times A$. On the other hand, $C: \mathbb{K} \to \mathbb{R}$ is the one-step cost function, and $P = [p_{xy}(\cdot)]$ is the controlled transition law. The interpretation of $M$ is as follows: At each time $t \in \mathbb{N}$ the state of a dynamical system is observed, say $X_t = x \in S$, and an action $A_t = a \in A(x)$ is chosen. Then a cost $C(x,a)$ is incurred and, regardless of the previous states and actions, the state of the system at time $t+1$ will be $X_{t+1} = y \in S$ with probability $p_{xy}(a)$; this is the Markov property of the decision model.

*Policies.* For each $t \in \mathbb{N}$ the space $\mathbb{H}_t$ of admissible histories up to time $t$ is recursively defined by $\mathbb{H}_0 := S$, and $\mathbb{H}_t := \mathbb{K} \times \mathbb{H}_{t-1}$ for $t \geq 1$. A generic element of $\mathbb{H}_t$ is denoted by $\mathbf{h}_t = (x_0, a_0, x_1, a_i, \ldots, x_{t-1}, a_{t-1}, x_t)$, where $x_n \in S$ for $n \leq t$, and $a_i \in A(x_i)$ for $i < t$. A policy $\pi = \{\pi_t\}$ is a special sequence of stochastic kernels: For each $t \in \mathbb{N}$ and $\mathbf{h}_t \in \mathbb{H}_t$, $\pi_t(\cdot | h_t)$ is a probability measure on $A$ concentrated on $A(x_t)$. The class of all policies is denoted by $\mathcal{P}$. Given the policy $\pi \in \mathcal{P}$ used to drive the system and the initial state $X_0 = x \in S$, the distribution of the state-action process $\{(X_t, A_t)\}$ is uniquely determined via Ionescu Tulcea's theorem (see, e.g., [[15]] or [[18]]); such a distribution will be represented by $P_x^\pi$, whereas $E_x^\pi$ stands for the corresponding expectation operator. Throughout the remainder $I_t$ denotes the information vector up to time $t$, which is given by

$$I_0 = X_0 \quad \text{and} \quad I_t := (X_0, A_0, \ldots, X_{t-1}, A_{t-1}, X_t), \qquad t = 1, 2, 3, \ldots.$$

Next, define $\mathbb{F} := \prod_{x \in S} A(x)$ so that $\mathbb{F}$ consists of all (choice) functions $f: S \to A$ satisfying that $f(x) \in A(x)$ for each $x \in S$. A policy $\pi$ is stationary if there exists $f \in \mathbb{F}$ such that, when the system evolves under $\pi$, at each time $t \in \mathbb{N}$ the action applied is determined by $A_t = f(X_t)$; the class of stationary policies is naturally identified with $\mathbb{F}$ and, with this convention, $\mathbb{F} \subset \mathcal{P}$.

*Performance index.* As already noted, the controller is assumed to be risk-averse with constant risk sensitivity $\lambda > 0$, that is, when facing a random cost $Y$, she grades it through $E[e^{\lambda Y}]$. The certain equivalent of the random variable $Y$ is the (possibly extended) real number defined by

$$E(\lambda, Y) := \frac{1}{\lambda} \log(E[e^{\lambda Y}]),$$

so that $e^{\lambda E(\lambda, Y)} = E[e^{\lambda Y}]$, and then the controller is indifferent between incurring the random cost $Y$ or paying the certain equivalent $E(\lambda, Y)$ for sure.

When the system evolves under $\pi \in \mathcal{P}$ and $x \in S$ is the initial state, $J_n(\lambda, \pi, x)$ denotes the certain equivalent of the total cost incurred before time $n > 0$, that is,

$$(2.1) \qquad J_n(\lambda, \pi, x) := \frac{1}{\lambda} \log \left( E_x^\pi \left[ \exp \left\{ \lambda \sum_{t=0}^{n-1} C(X_t, A_t) \right\} \right] \right),$$

whereas the (long-run expected $\lambda$-sensitive) average cost under $\pi$ starting at $x$ is defined by

$$(2.2) \qquad J(\lambda, \pi, x) := \limsup_{n \to \infty} \frac{1}{n} J_n(\lambda, \pi, x).$$

The optimal ($\lambda$-sensitive) average cost at state $x$ is given by

$$(2.3) \qquad J^*(\lambda, x) := \inf_\pi J(\lambda, \pi, x),$$

and a policy $\pi^* \in \mathcal{P}$ is optimal if $J(\lambda, \pi^*, x) = J^*(\lambda, x)$ for every $x \in S$. Given $\varepsilon > 0$, a policy $\pi$ is $\varepsilon$-optimal at state $x \in S$ if $J(\lambda, \pi, x) \leq J^*(\lambda, x) + \varepsilon$; if the policy $\pi$ is $\varepsilon$-optimal at every state, then $\pi$ is $\varepsilon$-optimal. The following simultaneous Doeblin condition will be assumed throughout the sequel.

ASSUMPTION 2.1. There exists a state $z \in S$ and $M \in (0, \infty)$ such that

$$E_x^f[T] \leq M, \qquad x \in S, f \in \mathbb{F},$$

where

$$(2.4) \qquad T := \min\{n > 0 | X_n = z\}$$

is the first positive arrival time to state $z$ and, by convention, the minimum of the empty set is $\infty$.

*The problem.* As already mentioned, the main objective of the paper is to provide a characterization of the optimal value function $J^*(\lambda, \cdot)$ for arbitrary $\lambda > 0$. This problem has recently received considerable attention in the literature and, under the above simultaneous Doeblin condition, the results already established can be described as follows: if $\lambda > 0$ is sufficiently small, then the optimal value function $J^*(\lambda, \cdot)$ is constant and, moreover, its value $\gamma$ is the unique real number for which there exists $h: S \to \mathbb{R}$ satisfying the optimality equation

$$(2.5) \qquad e^{\lambda[\gamma + h(x)]} = \min_{a \in A(x)} \left[ e^{\lambda C(x,a)} \sum_y p_{xy}(a) e^{\lambda h(y)} \right], \qquad x \in S;$$

see [3, 5, 14]. Also, modulo an additive constant, the relative value function $h(\cdot)$ in this equation satisfies that for each $x \in S$,

$$(2.6) \begin{aligned} h(x) &= \inf_{\pi \in \mathcal{P}} \frac{1}{\lambda} \log\left( E_x^\pi \left[ \exp\left\{ \lambda \sum_{t=0}^{T-1} [C(X_t, A_t) - \gamma] \right\} \right] \right) \\ &= \inf_{\pi \in \mathcal{P}} \frac{1}{\lambda} \log\left( E_x^\pi \left[ \exp\left\{ \lambda \sum_{t=0}^{T-1} [C(X_t, A_t) - J^*(\lambda, X_t)] \right\} \right] \right), \end{aligned}$$

where $T$ is the hitting time in (2.4). However, the situation is substantially different when $\lambda > 0$ is arbitrary in that (i) Assumption 2.1 does not generally imply that $J^*(\lambda, \cdot)$ is constant, (ii) the rightmost term in (2.6) may be $\infty$ and, moreover, (iii) even when the optimal value function takes on a single value $\gamma$, it is not necessarily determined by (2.5). This potentially complex behavior, which does not occur under Assumption 2.1 when the performance index is the risk-neutral average cost, is illustrated in the following example along the lines of Example 2.1 in [6]. In all, this example shows that, when the risk sensitivity coefficient is large enough, the behavior of the system at transient states, which may be occupied only at "early stages," has an important and definite influence on its performance, establishing a remarkable difference with the risk-neutral case.

EXAMPLE 2.2. Let $S = \{0, 1, 2\}$ and $A = \{0, 1\}$. The sets of admissible actions are given by $A(0) = A(2) = \{0\}$ and $A(1) = \{0, 1\} = A$, whereas the cost function always satisfies $C(x, a) = x$. Finally, for some $\rho \in (0, 1)$, the transition law is determined by

$$p_{00}(0) = 1, \qquad p_{22}(0) = \rho^2 = 1 - p_{20}(0)$$

and

$$p_{12}(1) = 1, \qquad p_{11}(0) = \rho = 1 - p_{10}(0).$$

In this context it is not difficult to see that Assumption 2.1 is satisfied with $z = 0$. Now, let $f$ be the stationary policy determined by $f(1) = 0$ so that, since $a = 0$ is the unique action available at the absorbing state 0 and $C(0,0) = 0$, it follows that $J^*(\lambda, 0) = J(\lambda, f, 0) = 0$. Assume now that

$$e^\lambda \rho > 1. \tag{2.7}$$

Using that 0 is the unique available action at state 2 and that when the system leaves state 2 it reaches $z = 0$, where a null cost is incurred forever, it follows that $J^*(\lambda, 2) = J(\lambda, f, 2)$, whereas for each positive integer $n$,

$$\begin{aligned}
E_2^f &\left[ \exp\left\{ \lambda \sum_{t=0}^{n-1} C(X_t, A_t) \right\} \right] \\
&= \sum_{k=1}^n E_2^f \left[ \exp\left\{ \lambda \sum_{t=0}^{k-1} C(X_t, A_t) \right\} I[T = k] \right] \\
&\quad + E_2^f \left[ \exp\left\{ \lambda \sum_{t=0}^{n-1} C(X_t, A_t) \right\} I[T > n] \right] \\
&= \sum_{k=1}^n e^{2k\lambda} (\rho^2)^{k-1} (1 - \rho^2) + e^{2n\lambda} (\rho^2)^n \\
&= (e^\lambda \rho)^{2n} + e^{2\lambda}(1 - \rho^2) \frac{(e^\lambda \rho)^{2n} - 1}{(e^\lambda \rho)^2 - 1}
\end{aligned}$$

and then (2.1) and (2.2) together lead to $J(\lambda, f, 2) = \frac{1}{\lambda} \log[(e^{2\lambda} \rho^2)] = \frac{2}{\lambda} \log(e^\lambda \rho) > 0$. A similar argument shows that $J(\lambda, f, 1) = \frac{1}{\lambda} \log(e^\lambda \rho) > 0$ and, since applying action 1 at state 1 produces a transition to state 2, where the optimal average cost is $\frac{2}{\lambda} \log(e^\lambda \rho) > J(\lambda, f, 1)$, it follows that $f$ is also optimal at state 1. In short, under (2.7),

$$J^*(\lambda, 0) = 0 < \frac{1}{\lambda} \log(e^\lambda \rho) = 1 + \frac{\log(\rho)}{\lambda} = J^*(\lambda, 1) < 2J^*(\lambda, 1) = J^*(\lambda, 2).$$

Notice that the system will be ultimately absorbed by state $z = 0$ but, when (2.7) holds, the costs incurred at the transient states have a definite influence on the performance of the system. Assume now that the initial state is $X_0 = 2$. From the specification of the model, it follows that $X_t = 2$ for $t < T$, with $T$ as in (2.4) with $z = 0$, and in this case $C(X_t, A_t) - J^*(\lambda, X_t) = 2 - 2(1 + \log(\rho)/\lambda) = -2\log(\rho)/\lambda$, so that $\lambda \sum_{t=0}^{T-1}[C(X_t, A_t) - J^*(\lambda, X_t)] = -2T\log(\rho)$. Therefore, the relative value function at state $x = 2$, given by the rightmost term in (2.6), is $h(2) = E_2^f[e^{-2T\log(\rho)}] = E_2^f[\rho^{-2T}] = \sum_{k=1}^\infty \rho^{-2k}(\rho^2)^k(1 - \rho^2) = \infty$; similarly, it can be established that $h(1) = \infty$.

On the other hand, it is interesting to observe that there is not any function $h: S \to \mathbb{R}$ satisfying that

$$(2.8) \qquad e^{\lambda J^*(\lambda,2)+\lambda h(2)} \geq e^{\lambda C(2,0)} \sum_y p_{2y}(0) e^{\lambda h(y)};$$

indeed, the left-hand side of this inequality is $e^{2\lambda}\rho^2 e^{\lambda h(2)}$, whereas the right-hand side satisfies $e^{\lambda C(2,0)}[p_{22}(0)e^{\lambda h(2)} + p_{21}(0)e^{\lambda h(0)}] > e^{\lambda C(2,0)} p_{22}(0) e^{\lambda h(2)} = e^{2\lambda}\rho^2 e^{\lambda h(2)}$. When the risk sensitivity coefficient satisfies $e^\lambda \rho = 1$, similar calculations yield that (i) $J^*(\lambda, \cdot) \equiv 0 = \gamma$, (ii) the relative value function $h$ in (2.6) is $\infty$ at $x = 1$ and 2, and (iii) inequality (2.8) is not satisfied by any function $h: S \to \mathbb{R}$; in particular, even in this case in which the optimal average cost is constant, the optimality equation (2.5) does not have a solution. Finally, if $\lambda$ satisfies that $e^\lambda \rho < 1$, which in this example is the precise meaning of "if $\lambda$ is sufficiently small," the optimal value function is identically $0 = \gamma$, the relative value function in (2.6) is finite, and the pair $(\gamma, h(\cdot))$ satisfies the optimality equation (2.5); see [3] or [14] for these latter assertions.

*The characterization theorem.* The main result of this work, which is formally stated as Theorem 3.5 in the following section, provides a characterization of $J^*(\lambda, \cdot)$ covering the diversity of possible behaviors illustrated in Example 2.2. For each $\lambda > 0$, this theorem determines the optimal value function in terms of a class of functions $\mathcal{G}$, and establishes that $J^*(\lambda, \cdot)$ is the infimum of such a family. This conclusion is similar to the characterization of the optimal (risk-neutral) total expected cost $V^*$ for MDPs with nonnegative cost function; in this latter case, $V^*$ is the infimum of all nonnegative functions $W$ defined on the state space and satisfying $W \geq \mathcal{D}W$, where $\mathcal{D}$ is the corresponding dynamic programming operator; see [18] for details. However, for the risk-sensitive average criterion in this work, the construction of family $\mathcal{G}$ involves *two conditions*, resembling the two equations that characterize the optimal risk-neutral average cost in multichain MDPs, for which the optimal performance index is not necessarily constant (see, e.g., Chapter 9 in [18]). The first restriction imposed on the members of $\mathcal{G}$ reflects a fundamental property of the risk-sensitive average index, namely, if the system is driven by a "good" policy, then $\{J^*(\lambda, X_t)\}$ is nonincreasing for almost all sample trajectories. This property is a consequence of Lemma 3.1 in the following section, establishing that the optimal value function satisfies a min–max equation, and the first condition imposed on the members of $\mathcal{G}$ is to satisfy such an equality. The second condition on a function $g \in \mathcal{G}$ is motivated by the optimality equation that, at least formally, is associated with this optimal control problem. This condition guarantees that $g$ is really an upper bound of $J^*(\lambda, \cdot)$; it was also used in [6] to analyze the uncontrolled case, and requires the existence of a (deviation) function $h: S \to \mathbb{R}$ such that the pair $(g(\cdot), h(\cdot))$ satisfies (3.4), which is analogous to the condition $W \geq \mathcal{D}W$ mentioned above.

*Outline of the argument.* As might be expected from the diversity illustrated in Example 2.2, characterizing $J^*(\lambda, \cdot)$ for arbitrary $\lambda > 0$ is a somewhat technical task, so that it is convenient to give a brief outline of the argument used to achieve this goal. In Section 3 the basic min–max equation satisfied by $J^*(\lambda, \cdot)$ is established, and then the family of functions $\mathcal{G}$ is introduced. Next, it is shown that the optimal $\lambda$-sensitive average cost is a lower bound of $\mathcal{G}$, and the characterization result of $J^*(\lambda, \cdot)$ as the infimum of $\mathcal{G}$ is stated as Theorem 3.5. As it will be noted below, in general $J^*(\lambda, \cdot)$ does not belong to $\mathcal{G}$, but the strategy to establish Theorem 3.5 consists in showing that, for each $\alpha \in (0, 1)$, the function $g(\cdot) = \alpha J^*(\lambda, \cdot) + (1 - \alpha)\|C\|$ lies in $\mathcal{G}$, from which Theorem 3.5 follows immediately. The main difficulty in establishing this inclusion is to prove that there exists a deviation function $h: S \to \mathbb{R}$ such that the second condition in the definition of family $\mathcal{G}$ is satisfied. In Definition 4.1 a candidate $h$ for the deviation function for the function $g$ above is introduced, and from that point onward, the effort is mainly dedicated to establishing that $h(\cdot)$ is a finite function, a fact that is proved in two steps: In Theorem 4.4 it is shown that $h$ is finite at the points $x$ where the optimal value function is minimized, whereas in Theorem 5.1 this conclusion is extended to the whole state space. The argument in this part relies heavily on the following property: Under an $\varepsilon$-optimal policy with $\varepsilon > 0$ small enough, along almost all trajectories the optimal value function is dominated by its value at the initial state. Section 6 concerns a last technical point on the function $h$ introduced in Definition 4.1, namely, that $h(z)$ is nonpositive, where $z$ is as in Assumption 2.1. After the preliminaries in Sections 4–6, Theorem 3.5 is finally proved in Section 7.

Before leaving this section, it is convenient to point out the following observation.

REMARK 2.3. (i) Given $\varepsilon > 0$, an $\varepsilon$-optimal policy exists. Indeed, from the definition of $J^*(\lambda, \cdot)$ in (2.3), it follows that for each $x \in S$ there exists a policy $\pi^x \in \mathcal{P}$ which is $\varepsilon$-optimal at $x$, that is,

$$J(\lambda, \pi^x, x) \leq J^*(\lambda, x) + \varepsilon$$

and a new policy $\pi$ can be defined as follows: For each $t \in \mathbb{N}$ and $\mathbf{h}_t \in \mathbb{H}_t$, $\pi_t(\cdot|\mathbf{h}_t) = \pi_t^{x_0}(\cdot|\mathbf{h}_t)$. A controller driving the system according to $\pi$ first determines the initial state, and then picks the actions according to $\pi^x$ if $X_0 = x$ is observed. From this construction it follows that the equality

$$E_x^\pi\left[\exp\left\{\lambda \sum_{t=0}^n C(X_t, A_t)\right\}\right] = E_x^{\pi^x}\left[\exp\left\{\lambda \sum_{t=0}^n C(X_t, A_t)\right\}\right]$$

is always valid, and then (2.1) and (2.2) together yield that $J(\lambda, \pi, x) = J(\lambda, \pi^x, x) \leq J^*(\lambda, x) + \varepsilon$ for every state $x$, so that $\pi$ is $\varepsilon$-optimal.

(ii) From (2.1)–(2.3) it is not difficult to see that $-\|C\| \leq J^*(\lambda, \cdot) \leq \|C\|$.

**3. Min–max equation and main result.** According to the program outlined above, in this section the characterization result for the optimal value function is stated. First, it is shown in the next lemma that the fundamental min–max equation is satisfied by the optimal value function, and such an equality is used as one of the requirements in the definition of the family of functions $\mathcal{G}$ involved in the characterization of $J^*(\lambda, \cdot)$.

LEMMA 3.1. *For each $\lambda > 0$, the function $J^*(\lambda, \cdot)$ in (2.3) satisfies the following min–max equation:*

$$J^*(x) = \min_{a \in A(x)} \max\{J^*(y) | p_{xy}(a) > 0\}, \qquad x \in S.$$

PROOF. Let $(x, a) \in \mathbb{K}$ and $\varepsilon > 0$ be arbitrary but fixed, and let $\pi \in \mathcal{P}$ be an $\varepsilon$-optimal policy (see Remark 2.3). Next, select a policy $f \in \mathbb{F}$ satisfying that $f(x) = a$, and define the new policy $\tilde{\pi} \in \mathcal{P}$ as follows: $\tilde{\pi}_0(\{f(x_0)\}|x_0) = 1$ for each $x_0 \in S$, whereas for each $t \in \mathbb{N}$ and $\mathbf{h}_{t+1} \in \mathbb{H}_{t+1}$,

$$\tilde{\pi}_{t+1}(\cdot|\mathbf{h}_{t+1}) = \pi_t(\cdot|x_1, a_1, \ldots, x_{t+1}).$$

When the system is driven by $\tilde{\pi}$, the action applied at time zero is selected using $f$, whereas from time 1 onwards, the controls are picked using the $\varepsilon$-optimal policy $\pi$ as if the decision process had started again at time 1. The Markov property and (2.1) together yield that for every positive integer $n$,

$$\begin{aligned} e^{\lambda J_{n+1}(\lambda, \tilde{\pi}, x)} &= E_x^{\tilde{\pi}}\left[\exp\left\{\lambda \sum_{t=0}^{n} C(X_t, A_t)\right\}\right] \\ &= e^{\lambda C(x, f(x))} \sum_y p_{xy}(f(x)) E_y^{\pi}\left[\exp\left\{\lambda \sum_{t=0}^{n-1} C(X_t, A_t)\right\}\right] \\ &= e^{\lambda C(x, a)} \sum_y p_{xy}(a) e^{\lambda J_n(\lambda, \pi, y)} \end{aligned}$$

so that

$$(3.1) \quad \frac{J_{n+1}(\lambda, \tilde{\pi}, x)}{n+1} \leq \frac{C(x, a)}{n+1} + \frac{n}{n+1} \log\left(\left[\sum_y p_{xy}(a) e^{\lambda J_n(\lambda, \pi, y)}\right]^{1/(\lambda n)}\right).$$

On the other hand, since $\pi$ is $\varepsilon$-optimal and $S$ is finite, it follows that for some $n_0 \in \mathbb{N}$,

$$J_n(\lambda, \pi, \cdot) \leq n(J^*(\lambda, \cdot) + \varepsilon), \qquad n \geq n_0.$$

Therefore, $\sum_y p_{xy}(a)e^{\lambda J_n(\lambda,\pi,y)} \leq \sum_y p_{xy}(a)e^{\lambda n(J^*(\lambda,y)+\varepsilon)}$ when $n \geq n_0$, and it follows that

$$\limsup_{n\to\infty}\left[\sum_y p_{xy}(a)e^{\lambda J_n(\lambda,\pi,y)}\right]^{1/(\lambda n)} \leq \limsup_{n\to\infty}\left[\sum_y p_{xy}(a)e^{\lambda n(J^*(\lambda,y)+\varepsilon)}\right]^{1/(\lambda n)}$$
$$= \max\{e^{J^*(\lambda,y)+\varepsilon}|p_{xy}(a)>0\}$$
$$= e^{\max\{J^*(\lambda,y)+\varepsilon|p_{xy}(a)>0\}},$$

where the second equality is due to the fact that the exponential function is increasing. Combining this with (2.2), after taking limit superior as $n$ goes to $\infty$ in (3.1) it follows that

$$J^*(\lambda,x) \leq J(\lambda,\tilde\pi,x) \leq \max\{J^*(\lambda,y)+\varepsilon|p_{xy}(a)>0\};$$

see (2.3) for the first inequality. Recalling that $\varepsilon > 0$ is arbitrary, this yields

$$J^*(\lambda,x) \leq \max\{J^*(\lambda,y)|p_{xy}(a)>0\},$$

a relation that, since $(x,a) \in \mathbb{K}$ is arbitrary, implies

$$(3.2) \qquad J^*(\lambda,x) \leq \min_{a \in A(x)} \max\{J^*(\lambda,y)|p_{xy}(a)>0\}, \qquad x \in S.$$

To establish the reverse inequality let $x \in S$ and $\pi \in \mathcal{P}$ be arbitrary. Select $b \in A(x)$ satisfying $\pi_0(\{b\}|x) > 0$, and let $y \in S$ be such that $p_{xy}(b) > 0$. Combining (2.1) with the Markov property, it follows that for every positive integer $n$,

$$e^{\lambda J_{n+1}(\lambda,\pi,x)} = E_x^\pi\left[\exp\left\{\lambda \sum_{t=0}^n C(X_t,A_t)\right\}\right]$$
$$\geq E_x^\pi\left[\exp\left\{\lambda \sum_{t=0}^n C(X_t,A_t)\right\}I[A_0=b,X_1=y]\right]$$
$$= \pi_0(\{b\}|x)p_{xy}(b)e^{\lambda C(x,b)}E_y^\delta\left[\exp\left\{\lambda \sum_{t=0}^{n-1} C(X_t,A_t)\right\}\right]$$
$$= \pi_0(\{b\}|x)p_{xy}(b)e^{\lambda C(x,b)}e^{\lambda J_n(\lambda,\delta,y)},$$

where the "shifted" policy $\delta$ is defined as follows: For every $t \in \mathbb{N}$ and $\mathbf{h}_t \in \mathbb{H}_t$, $\delta_t(\cdot|\mathbf{h}_t) = \pi_{t+1}(\cdot|x,b,\mathbf{h}_t)$. Therefore,

$$\frac{J_{n+1}(\lambda,\pi,x)}{n+1} \geq \frac{1}{\lambda(n+1)}\log(\pi_0(\{b\}|x)p_{xy}(b)e^{\lambda C(x,b)}) + \frac{n}{n+1}\frac{J_n(\lambda,\delta,y)}{n},$$

and taking limit superior as $n$ goes to $\infty$, it follows that

$$J(\lambda,\pi,x) \geq J(\lambda,\delta,y) \geq J^*(\lambda,y);$$

see (2.2) and (2.3). Since the state $y$ satisfying $p_{xy}(b) > 0$ is arbitrary, this implies that

$$J(\lambda, \pi, x) \geq \max\{J^*(\lambda, y)|p_{xy}(b) > 0\},$$

and then

$$J(\lambda, \pi, x) \geq \min_{a \in A(x)} \max\{J^*(\lambda, y)|p_{xy}(a) > 0\}.$$

Since this holds for every $\pi \in \mathcal{P}$ and $x \in S$, (2.2) yields that

$$J^*(\lambda, x) \geq \min_{a \in A(x)} \max\{J^*(\lambda, y)|p_{xy}(a) > 0\}, \qquad x \in S,$$

and the result follows combining this inequality with (3.2). $\square$

DEFINITION 3.2. The class $\mathcal{G}$ consists of all functions $g \in \mathcal{B}(S)$ satisfying the following conditions:

(i) For each $x \in S$

(3.3) $$g(x) = \min_{a \in A(x)} \max\{g(y)|p_{xy}(a) > 0\}.$$

(ii) There exists a function $h \in \mathcal{B}(S)$, possibly depending on $g$, such that

(3.4) $$e^{\lambda g(x) + \lambda h(x)} \geq \min_{a \in B_g(x)} \left[e^{\lambda C(x,a)} \sum_y p_{xy}(a) e^{\lambda h(y)}\right], \qquad x \in S,$$

where

(3.5) $$B_g(x) := \{a \in A(x)|g(x) = \max\{g(y)|p_{xy}(a) > 0\}\};$$

a function $h(\cdot)$ satisfying (3.4) will be referred to as a deviation function associated to $g(\cdot)$.

REMARK 3.3. (i) Given $g \in \mathcal{B}(S)$ satisfying (3.3), the finiteness of the action sets $A(x)$ ensures that each set $B_g(x)$ is nonempty.

(ii) Family $\mathcal{G}$ is nonempty. In fact, if $g(\cdot) = \|C\|$, then $g \in \mathcal{G}$, since (3.3) is clearly satisfied by this function, whereas (3.4) holds with $h(\cdot) = 0$.

The following lemma shows that the optimal value function is dominated by each member of $\mathcal{G}$.

LEMMA 3.4. (i) *Suppose that $g: S \to \mathbb{R}$ satisfies (3.3) and for each $x \in S$ let $B_g(x) \subset A(x)$ be as in (3.5). Given $x \in S$, assume that the policy $\delta \in \mathcal{P}$ satisfies that $P_x^\delta[A_r \in B_g(X_r)] = 1$ for every $r \in \mathbb{N}$. In this case, when $x$ is the initial state and the system is driven by $\delta$, the process $\{g(X_t)\}$ is nonincreasing almost surely. More precisely, for each $n \in \mathbb{N}$,*

$$g(X_{n+1}) \leq g(X_n) \leq \cdots \leq g(X_0) = g(x), \qquad P_x^\delta\text{-a.s.}$$

*Consequently,*

(ii) *Every $g \in \mathcal{G}$ is an upper bound of the optimal value function $J^*(\lambda, \cdot)$.*

PROOF. (i) Let $t \in \mathbb{N}$ be fixed, and suppose that $w, y \in S$ satisfy $P_x^\delta[X_t = w, X_{t+1} = y] > 0$. In this case there exists $a \in B_g(w)$ such that $P_x^\delta[X_t = w, A_t = a, X_{t+1} = y] > 0$, since $P_x^\delta[A_t \in B_g(X_t)] = 1$, and then

$$0 < P_x^\delta[X_t = w, A_t = a, X_{t+1} = y]$$
$$= P_x^\delta[X_{t+1} = y | X_t = w, A_t = a] P_x^\delta[X_t = w, A_t = a]$$
$$= p_{wy}(a) P_x^\delta[X_t = w, A_t = a],$$

where the second equality is due to the Markov property; therefore,

$$p_{wy}(a) > 0.$$

On the other hand, from (3.5), the inclusion $a \in B_g(w)$ yields that

$$g(w) = \max\{g(z) | p_{wz}(a) > 0\}$$

and combining this with the above inequality, it follows that $g(w) \geq g(y)$. In short, it has been shown that

$$P_x^\delta[X_t = w, X_{t+1} = y] > 0 \implies g(w) \geq g(y)$$

and it follows that $P_w^\delta[g(X_{t+1}) \leq g(X_t)] = 1$. Since $t \in \mathbb{N}$ is arbitrary, this yields that $P_w^\delta[g(X_{n+1}) \leq g(X_n) \leq \cdots \leq g(X_0)] = 1$ for each $n \in \mathbb{N}$.

(ii) Let $g \in \mathcal{G}$ be arbitrary, and select $h \in \mathcal{B}(S)$ as in (3.4). For each $x \in S$, let $f(x) \in B_g(x)$ be a minimizer of the term within brackets in (3.4), so that for every $x \in S$,

$$e^{\lambda g(x) + \lambda h(x)} \geq e^{\lambda C(x, f(x))} \sum_y p_{xy}(f(x)) e^{\lambda h(y)},$$

which is equivalent to $e^{\lambda h(x)} \geq E_x^f[e^{\lambda[C(X_0, A_0) - g(X_0)]} e^{\lambda h(X_1)}]$; from this point, an induction argument yields that

$$e^{\lambda h(x)} \geq E_x^f\left[\exp\left\{\lambda \sum_{t=0}^n [C(X_t, A_t) - g(X_t)]\right\} e^{\lambda h(X_{n+1})}\right], \qquad x \in S, n \in \mathbb{N}.$$

Observe that, by part (i), under the action of policy $f$ the inequalities

$$g(X_n) \leq g(X_{n-1}) \leq \cdots \leq g(X_1) \leq g(X_0)$$

hold with probability 1 regardless of the initial state. Therefore, $\sum_{t=0}^n[C(X_t, A_t) - g(X_t)] \geq \sum_{t=0}^n C(X_t, A_t) - (n+1)g(X_0)$, so that for each $n \in \mathbb{N}$ and $x \in S$,

$$e^{\lambda h(x)} \geq E_x^f\left[\exp\left\{\lambda \sum_{t=0}^n C(X_t, A_t) - (n+1)g(X_0)\right\} e^{\lambda h(X_{n+1})}\right]$$

$$= e^{-\lambda(n+1)g(x)} E_x^f \left[ \exp\left\{ \lambda \sum_{t=0}^n C(X_t, A_t) \right\} e^{\lambda h(X_{n+1})} \right]$$

$$\geq e^{-\lambda(n+1)g(x) - \lambda\|h\|} E_x^f \left[ \exp\left\{ \lambda \sum_{t=0}^n C(X_t, A_t) \right\} \right]$$

$$\geq e^{-\lambda(n+1)g(x) - \lambda\|h\|} e^{\lambda J_{n+1}(\lambda, f, x)}.$$

Hence,
$$g(x) + \frac{\|h\| + h(x)}{n+1} \geq \frac{J_{n+1}(\lambda, f, x)}{n+1},$$

and taking limit superior as $n$ goes to $\infty$, this yields that $g(x) \geq J(\lambda, f, x) \geq J^*(\lambda, x)$; see (2.2) and (2.3). Since $x \in S$ is arbitrary, it follows that $g(\cdot) \geq J^*(\lambda, \cdot)$. $\square$

According to this result, the functional $J^*(\lambda, \cdot)$ is a lower bound for each member of $\mathcal{G}$. On the other hand, although $J^*(\lambda, \cdot)$ satisfies (3.3), by Lemma 3.1, in general this optimal value function does not belong to $\mathcal{G}$. Indeed, in the context of Example 2.2, it was shown that when $e^\lambda \rho \geq 1$ there is not any function $h$ such that (2.8) is satisfied, and this implies that the second part of Definition 3.2 fails for the function $J^*(\lambda, \cdot)$. However, under Assumption 2.1, the main result of this work asserts that $J^*(\lambda, \cdot)$ is the largest lower bound of $\mathcal{G}$.

THEOREM 3.5.  *Under Assumption* 2.1, *for each* $x \in S$,
$$J^*(\lambda, x) = \inf_{g \in \mathcal{G}} g(x).$$

This result extends Theorem 2.2 in [6] where the uncontrolled case was analyzed. The somewhat technical proof of this theorem will be given in Section 7 after establishing the necessary technical preliminaries in the following three sections. Essentially, although it cannot be ensured that the optimal value function is a member of $\mathcal{G}$, the idea is to prove that, for each $\alpha \in (0, 1)$, the function $g$ specified by

(3.6) $$g(\cdot) := \alpha J^*(\lambda, \cdot) + (1 - \alpha)\|C\|$$

lies in $\mathcal{G}$, a fact that immediately yields Theorem 3.5. Using Lemma 3.1 it is not difficult to see that this function $g$ satisfies the min–max equation (3.3) and then, to establish the inclusion $g \in \mathcal{G}$ it is sufficient to show that there exists a deviation function $h \in \mathcal{B}(S)$ associated to $g$, so that the pair $(g, h)$ satisfies (3.4). The proof of this existence result requires an important technical effort that is presented in the following three sections. Throughout the remainder Assumption 2.1 is supposed to hold even without explicit reference, and $\alpha \in (0, 1)$ is arbitrary but fixed.

**4. Deviation function.** In this section a candidate $h(\cdot)$ for the deviation function of the function $g$ in (3.6) is introduced and, as already mentioned, a major objective is to show that such a function is finite. Although this goal is finally achieved later, the main result of this section, stated as Theorem 4.4, is a first step in this direction.

DEFINITION 4.1. (i) For each $x \in S$, define $B^*(x) := B_{J^*(\lambda,\cdot)}(x)$; see (3.5).

(ii) The class $\mathcal{P}^*$ consists of all policies $\pi \in \mathcal{P}$ satisfying
$$P_x^\pi[A_t \in B^*(X_t)] = 1, \qquad x \in S, t \in \mathbb{N}.$$

(iii) Given a fixed real number $\alpha \in (0,1)$, the corresponding deviation function $h: S \to [-\infty, \infty]$ is defined as

$$h(x) = \inf_{\pi \in \mathcal{P}^*} \frac{1}{\lambda} \log \left( E_x^\pi \left[ \exp \left\{ \lambda \alpha \sum_{t=0}^{T-1} [C(X_t, A_t) - J^*(\lambda, X_t)] \right\} \right] \right),$$
(4.1)
$$x \in S,$$

where $T$ is the first positive passage time to state $z$; see (2.4).

Notice that $a \in B^*(x)$ if and only if $J^*(\lambda, x) = \max\{J^*(\lambda, y) | p_{xy}(a) > 0\}$, and that $\mathcal{P}^*$ is the class of policies for the MDP $\langle S, A, \{B^*(x)\}, C, P \rangle$, which is obtained by restricting the set of admissible actions at state $x$ to the subset $B^*(x)$. On the other hand, observe that the factor $\lambda \alpha$ is used in the exponential inside the expectation in (4.1); when $\alpha = 1$, it is not difficult to see from Example 2.2 that $h(\cdot)$ may take on an infinite value at some points. However, in the present case in which $\alpha$ lies in $(0,1)$, it will be proved that $h(\cdot)$ is finite. The key tool in the argument leading to this goal is the following consequence of Assumption 2.1.

LEMMA 4.2. *Under Assumption 2.1 there exist $\beta \in (0,1)$ and $\beta_0 > 0$ such that*
$$P_x^\pi[T \geq n] \leq \beta_0 \beta^n, \qquad x \in S, \pi \in \mathcal{P}, n \in \mathbb{N}.$$

A proof of this lemma can be seen, for instance, in [16] or [19]. Using this result, it is now shown that $-\infty$ is not a value of the function $h(\cdot)$ in (4.1).

LEMMA 4.3. *Let $\alpha \in (0,1)$ be fixed and let $\beta_0$ and $\beta$ be as in Lemma 4.2.*

(i) *There exists a positive constant $B_0$ such that for each $x \in S$ and $\pi \in \mathcal{P}$,*
$$E_x^\pi \left[ \exp \left\{ \lambda \alpha \sum_{t=0}^{T-1} [C(X_t, A_t) - J^*(\lambda, X_t)] \right\} \right] \geq B_0.$$

*Consequently:*

(ii) $h(\cdot) > -\infty$.

(iii) *Set*

$$\xi_0 = -\frac{(1-\alpha)\log(\beta)}{\lambda\alpha}. \tag{4.2}$$

*If $\pi \in \mathcal{P}$ is $\varepsilon$-optimal, where $\varepsilon \in (0, \xi_0)$, and $x \in S$ is such that*

$$J^*(\lambda, x) = \gamma,$$

*then*

$$E_x^\pi\left[\exp\left\{\lambda\alpha\sum_{t=0}^{T-1}[C(X_t, A_t) - \gamma]\right\}\right] < \infty.$$

PROOF. (i) Let $N_0 \in \mathbb{N}$ be such that $\beta_0\beta^{N_0+1} < 1/2$, and observe that the inequality

$$P_x^\pi[T \leq N_0] \geq \tfrac{1}{2}$$

always holds by Lemma 4.2. On the other hand, using Remark 2.3(ii), it follows that

$$\sum_{t=0}^{T-1}[C(X_t, A_t) - J^*(\lambda, X_t)] \geq -2\|C\|T,$$

so that for every $x \in S$ and $\pi \in \mathcal{P}$,

$$E_x^\pi\left[\exp\left\{\lambda\alpha\sum_{t=0}^{T-1}[C(X_t, A_t) - J^*(\lambda, X_t)]\right\}\right]$$

$$\geq E_x^\pi[e^{-2\lambda\alpha\|C\|T}]$$

$$\geq E_x^\pi[e^{-2\lambda\|C\|T}] \geq \sum_{k=0}^{N_0} P_x^\pi[T = k]e^{-2\lambda\|C\|k}$$

and then

$$E_x^\pi\left[\exp\left\{\lambda\alpha\sum_{t=0}^{T-1}[C(X_t, A_t) - J^*(\lambda, X_t)]\right\}\right] \geq e^{-2\lambda\|C\|N_0}P_x^\pi[T \leq N_0]$$

$$\geq \frac{e^{-2\lambda\|C\|N_0}}{2} =: B_0.$$

(ii) Combining part (i) with (4.1), it follows that $h(\cdot) \geq -\log(B_0) > -\infty$.

(iii) Observe that Hölder's inequality implies

$$E_x^\pi\left[\exp\left\{\lambda\alpha\sum_{t=0}^{T-1}[C(X_t,A_t)-\gamma]\right\}\right]$$
$$=\sum_{n=1}^\infty E_x^\pi\left[\exp\left\{\lambda\alpha\sum_{t=0}^{n-1}[C(X_t,A_t)-\gamma]\right\}I[T=n]\right]$$
$$\leq\sum_{n=1}^\infty\left(E_x^\pi\left[\exp\left\{\lambda\sum_{t=0}^{n-1}[C(X_t,A_t)-\gamma]\right\}\right]\right)^\alpha(P_x^\pi[T=n])^{1-\alpha}$$

and then (2.1) and Lemma 4.2 together yield

$$E_x^\pi\left[\exp\left\{\lambda\alpha\sum_{t=0}^{T-1}[C(X_t,A_t)-\gamma]\right\}\right]\leq\beta_0^{(1-\alpha)}\sum_{n=1}^\infty e^{\lambda\alpha[J_n(\lambda,\pi,x)-n\gamma]}\beta^{n(1-\alpha)}.$$

Since $J^*(\lambda,x)=\gamma$ and $\pi$ is $\varepsilon$-optimal, it follows that $J_n(\lambda,\pi,x)/n\leq\gamma+\varepsilon$ when the positive integer $n$ is large enough, say $n>n_0$. Therefore,

$$e^{\lambda\alpha[J_n(\lambda,\pi,x)-n\gamma]}\beta^{n(1-\alpha)}\leq e^{\lambda\alpha\varepsilon n}\beta^{n(1-\alpha)}\leq(e^{\lambda\alpha\varepsilon}\beta^{(1-\alpha)})^n,\qquad n>n_0.$$

On the other hand, since $0<\varepsilon<\xi_0$, (4.2) implies that $e^{\lambda\alpha\varepsilon}\beta^{(1-\alpha)}<1$, so that the last two displayed relations together yield that $E_x^\pi[\exp\{\lambda\alpha\sum_{t=0}^{T-1}[C(X_t,A_t)-\gamma]\}]<\infty$. □

In contrast with the above argument used to establish the inequality $h(\cdot)>-\infty$, the proof of the inequality $h(\cdot)<\infty$ is substantially more technical. As a starting point, the main result of this section, stated in the following theorem, establishes that function $h(\cdot)$ is finite at the points where the optimal value function attains its minimum value.

THEOREM 4.4. *Let $\gamma_0$ be the minimum value of $J^*(\lambda,\cdot)$. In this case:*

(i) $J^*(\lambda,z)=\gamma_0$, *where $z$ is as in Assumption* 2.1.
(ii) *If $J^*(\lambda,x)=\gamma_0$, then $h(x)$ is finite.*

The proof of these results relies on the technical preliminaries in the following two lemmas; the first one provides a bound for $\{J^*(\lambda,X_t)\}$ when the system is driven by an $\varepsilon$-optimal policy.

LEMMA 4.5. *Let $\pi\in\mathcal{P}$, $x\in S$ and $r\in\mathbb{N}$ be arbitrary but fixed, and suppose that the vector $\tilde{\mathbf{h}}_r=(\tilde{x}_0,\tilde{a}_0,\ldots,\tilde{x}_{r-1},\tilde{a}_{r-1},\tilde{x}_r)\in\mathbb{H}_r$ satisfies*

$$P_x^\pi[I_r=\tilde{\mathbf{h}}_r]>0.$$

*In this case:*

(i) $J(\lambda, \pi, x) \geq J(\lambda, \delta, \tilde{x}_r) \geq J^*(\lambda, \tilde{x}_r)$, where the shifted policy $\delta$ is given by

(4.3) $\quad \delta_t(\cdot|\mathbf{h}_t) = \pi_{t+r}(\cdot|\tilde{x}_0, \tilde{a}_0, \ldots, \tilde{x}_{r-1}, \tilde{a}_{r-1}, \mathbf{h}_t), \quad t \in \mathbb{N}, \mathbf{h}_t \in \mathbb{H}_t.$

*Consequently:*

(ii) *If $\pi$ is $\varepsilon$-optimal at $x$, then for each $m \in \mathbb{N}$,*

$$J^*(\lambda, X_m) \leq J^*(\lambda, x) + \varepsilon, \qquad P_x^\pi\text{-a.s.}$$

PROOF. (i) Given an integer $n > r$, observe that

(4.4)
$$E_x^\pi\left[\exp\left\{\lambda \sum_{t=0}^n C(X_t, A_t)\right\}\right]$$
$$\geq e^{-\lambda r \|C\|} E_x^\pi\left[\exp\left\{\lambda \sum_{t=r}^n C(X_t, A_t)\right\} I[I_r = \tilde{\mathbf{h}}_r]\right].$$

On the other hand, an application of the Markov property yields that

$$E_x^\pi\left[\exp\left\{\lambda \sum_{t=r}^n C(X_t, A_t)\right\} I[I_r = \tilde{\mathbf{h}}_r] | I_r\right]$$
$$= I[I_r = \tilde{\mathbf{h}}_r] E_{\tilde{x}_r}^\delta\left[\exp\left\{\lambda \sum_{t=0}^{n-r} C(X_t, A_t)\right\}\right],$$

where policy $\delta$ is as in (4.3). Taking expectation with respect to $P_x^\pi$ in both sides of this equality, it follows that

$$E_x^\pi\left[\exp\left\{\lambda \sum_{t=r}^n C(X_t, A_t)\right\} I[I_r = \tilde{\mathbf{h}}_r]\right]$$
$$= P_x^\pi[I_r = \tilde{\mathbf{h}}_r] E_{\tilde{x}_r}^\delta\left[\exp\left\{\lambda \sum_{t=0}^{n-r} C(X_t, A_t)\right\}\right],$$

which combined with (4.4) leads to

$$E_x^\pi\left[\exp\left\{\lambda \sum_{t=0}^n C(X_t, A_t)\right\}\right]$$
$$\geq e^{-\lambda r \|C\|} P_x^\pi[I_r = \tilde{\mathbf{h}}_r] E_{\tilde{x}_r}^\delta\left[\exp\left\{\lambda \sum_{t=0}^{n-r} C(X_t, A_t)\right\}\right].$$

This inequality and (2.1) together imply that

$$\frac{J_{n+1}(\lambda, \pi, x)}{n+1} \geq \frac{\log(e^{-\lambda r \|C\|} P_x^\pi[I_r = \tilde{\mathbf{h}}_r])}{\lambda(n+1)} + \frac{n-r+1}{n+1} \frac{J_{n-r+1}(\lambda, \delta, \tilde{x}_r)}{n-r+1}$$

and taking limit superior as $n$ increases to $\infty$, this yields $J(\lambda, \pi, x) \geq J(\lambda, \delta, x_r) \geq J^*(\lambda, \tilde{x}_r)$; see (2.2).

(ii) Let $\pi \in \mathcal{P}$ be $\varepsilon$-optimal at $x$, and suppose that $P_x^\pi[X_m = y] > 0$. Observing that
$$[X_m = y] = \bigcup_{\mathbf{h}_m \in \mathbb{H}_m, x_m = y} [I_m = \mathbf{h}_m],$$
the finiteness of $\mathbb{H}_m$ implies that $P_x^\pi[I_m = \mathbf{h}_m] > 0$ for some $\mathbf{h}_m \in \mathbb{H}_m$ satisfying $x_m = y$. In this case, part (i) yields that $J^*(\lambda, y) \leq J(\lambda, \pi, x)$, so that $J^*(\lambda, y) \leq J^*(\lambda, x) + \varepsilon$, since $\pi$ is $\varepsilon$-optimal at $x$. In short
$$P_x^\pi[X_m = y] > 0 \implies J^*(\lambda, y) \leq J^*(\lambda, x) + \varepsilon,$$
and then $P_x^\pi[J^*(\lambda, X_m) \leq J^*(\lambda, x) + \varepsilon] = 1$. □

In the following lemma it is shown that, if $\varepsilon > 0$ is small enough, the set of minimizers of $J^*(\lambda, \cdot)$ is closed under the action of an $\varepsilon$-optimal policy and that, "essentially," such a policy belongs to the class $\mathcal{P}^*$ in Definition 4.1. The precise statement of these facts involves the following notation.

DEFINITION 4.6. (i) Define the positive number $\xi_1$ as follows:

(a) If $J^*(\lambda, \cdot)$ is constant, set $\xi_1 := 1$.
(b) If $J^*(\lambda, \cdot)$ is not constant, let $\gamma_i$, $i = 0, 1, \ldots, d$, be the different values of $J^*(\lambda, \cdot)$ arranged in increasing order:

(4.5) $$\gamma_0 < \gamma_1 < \cdots < \gamma_d.$$

In this case set
$$\xi_1 := \min\{\gamma_i - \gamma_{i-1} | i = 1, \ldots, d\}.$$

(ii) The positive number $\xi$ is given by
$$\xi = \min\{\xi_0, \xi_1\};$$
see (4.2).

REMARK 4.7. Observe that $J^*(\lambda, y) > J^*(\lambda, x)$ implies that $J^*(\lambda, y) \geq J^*(\lambda, x) + \xi_1$. Therefore,

if $0 < \varepsilon < \xi(\leq \xi_1)$, $\quad J^*(\lambda, x) + \varepsilon \geq J^*(\lambda, y) \implies J^*(\lambda, x) \geq J^*(\lambda, y)$.

LEMMA 4.8. Let $x \in S$ be such that $J^*(\lambda, x) = \gamma_0 = \min_y J^*(\lambda, y)$, and suppose that $\pi \in \mathcal{P}$ is $\varepsilon$-optimal at $x$, where $\varepsilon \in (0, \xi)$. In this case, for each $r \in \mathbb{N}$:

(i) $P_x^\pi[J^*(\lambda, X_r) = \gamma_0] = 1$.

(ii) $P_x^\pi[A_r \in B^*(X_r)] = 1$.

*Moreover:*

(iii) *there exists a policy $\delta \in \mathcal{P}^*$ such that, when the initial state is $x$, the distribution of the state-action process $\{(X_t, A_t)\}$ coincides under $\pi$ and $\delta$, that is, $P_x^\pi = P_x^\delta$.*

PROOF. (i) By Lemma 4.5, $P_x^\pi[J^*(\lambda, X_r) \leq J^*(\lambda, x) + \varepsilon] = 1$, whereas the inclusion $\varepsilon \in (0, \xi)$ yields that $[J^*(\lambda, X_r) \leq J^*(\lambda, x) + \varepsilon] \subset [J^*(\lambda, X_r) \leq J^*(\lambda, x)]$, by Remark 4.7, so that

$$P_x^\pi[J^*(\lambda, X_r) \leq J^*(\lambda, x)] = 1.$$

Since $J^*(\lambda, x) = \gamma_0$ is the minimum value of $J^*(\lambda, \cdot)$, it follows that $P_x^\pi[J^*(\lambda, X_r) = \gamma_0] = 1$.

(ii) Suppose that $P_x^\pi[A_r = a, X_r = w] > 0$. If $p_{wy}(a) > 0$, the Markov property yields that $P_x^\pi[X_{r+1} = y | X_r = w, A_r = a] = p_{wy}(a) > 0$, so that

$$0 < P_x^\pi[A_r = a, X_r = w] P_x^\pi[X_{r+1} = y | X_r = w, A_r = a]$$
$$= P_x^\pi[X_{r+1} = y, X_r = w, A_r = a]$$
$$\leq P_x^\pi[X_r = w, X_{r+1} = y]$$

and then part (i) yields that $J^*(\lambda, w) = J^*(\lambda, y) = \gamma_0$; since $y \in S$ satisfying $p_{wy}(a) > 0$, is arbitrary, it follows that $J^*(\lambda, w) = \max\{J^*(\lambda, y) | p_{xy}(a) > 0\}$, so that $a \in B^*(w)$; see Definition 4.1. Thus, $P_x^\pi[A_r = a, X_r = w] > 0 \Longrightarrow a \in B^*(w)$, so that

$$1 = \sum_{(w,a) \in \mathbb{K}} P_x^\pi[A_r = a, X_r = w]$$
$$= \sum_{(w,a) \in \mathbb{K}, a \in B^*(w)} P_x^\pi[A_r = a, X_r = w] = P_x^\pi[A_r \in B^*(X_r)].$$

(iii) Take a fixed stationary policy $f$ satisfying

$$f(y) \in B^*(y), \qquad y \in S,$$

and let the policy $\delta$ be determined as follows: For each $t \in \mathbb{N}$ and $\mathbf{h}_t \in \mathbb{H}_t$,

(4.6) $\delta_t(\cdot|\mathbf{h}_t) := \pi_t(\cdot|\mathbf{h}_t)$ if $x_0 = x, \delta_t(f(x_t)|\mathbf{h}_t) := 1$ when $x_0 \neq x$.

In this case, $\pi$ and $\delta$ coincide along trajectories starting at $x$, so that $P_x^\pi = P_x^\delta$, and then $P_x^\delta[A_t \in B^*(X_t)] = P_x^\pi[A_t \in B^*(X_t)] = 1$ for each $t \in \mathbb{N}$. Moreover, by the choice of $f$, $P_w^\delta[A_t \in B^*(X_t)] = 1$ always holds when $w \neq x$, and it follows that $\delta \in \mathcal{P}^*$; see Definition 4.1. $\square$

After the above preliminaries, the proof of the main result of this section can be established as follows.

PROOF OF THEOREM 4.4. Let $x \in S$ be a minimizer of $J^*(\lambda, \cdot)$, so that $J^*(\lambda, x) = \gamma_0$.

(i) By Lemma 4.2, there exists a positive integer $r$ such that $P_x^\pi[X_r = z] \geq P_x^\pi[T = r] > 0$, and then Lemma 4.8 (i) yields that $J^*(\lambda, z) = \gamma_0$.

(ii) Let $\pi$ be an $\varepsilon$-optimal policy, where $\varepsilon < \xi(\leq \xi_0)$; see (4.2) and Definition 4.6. In this case, using that $J^*(\lambda, x) = \gamma_0$, Lemma 4.3(iii) yields that $E_x^\pi[\exp\{\lambda\alpha \sum_{t=0}^{T-1}[C(X_t, A_t) - \gamma_0]\}] < \infty$; since $P_x^\pi[J^*(\lambda, X_r) = \gamma_0] = 1$ holds for every $r \in \mathbb{N}$, by Lemma 4.8(i), it follows that

$$E_x^\pi\left[\exp\left\{\lambda\alpha \sum_{t=0}^{T-1}[C(X_t, A_t) - J^*(\lambda, X_t)]\right\}\right] < \infty.$$

Now, using part (iii) in Lemma 4.8, select $\delta \in \mathcal{P}^*$ such that $P_x^\delta = P_x^\pi$, so that the above inequality yields

$$E_x^\delta\left[\exp\left\{\lambda\alpha \sum_{t=0}^{T-1}[C(X_t, A_t) - J^*(\lambda, X_t)]\right\}\right] < \infty$$

which, via Definition 4.1(iii) implies that $h(x) < \infty$; since $h(\cdot) > -\infty$, by Lemma 4.3(ii), it follows that $h(x)$ is finite. □

**5. Finiteness of the deviation function on the state space.** Following the program outlined in Section 2, the objective of this section is to extend the finiteness result in Theorem 4.4(ii) to the whole state space.

THEOREM 5.1. *For every $x \in S$, $h(x)$ is finite; see (4.1).*

Since $h(\cdot) > -\infty$, by Lemma 4.3(ii), to establish this result it is sufficient to show that $h(x) < \infty$ for every state $x \in S$. This latter inequality holds when $x$ is a minimizer of the optimal value function $J^*(\lambda, \cdot)$, by Theorem 4.4(ii), so that the deviation function is certainly finite when $J^*(\lambda, \cdot)$ is constant. Thus, to prove Theorem 5.1 it must be shown that $h(\cdot) < \infty$ when the optimal value function is not constant, and throughout the remainder of the section it is supposed that $J^*(\lambda, \cdot)$ assumes values $\gamma_i$, $i = 0, 1, \ldots, d$, where $d \geq 1$, which are arranged in increasing order; see (4.5). With this in mind, let the level set $G_i$ be given by

(5.1) $\qquad G_i := \{x \in S | J^*(\lambda, x) = \gamma_i\}, \qquad i = 0, \ldots, d.$

Notice that

(5.2) $$S = \bigcup_{i=0}^{d} G_i,$$

and define the exit time of set $G_i$ by

(5.3) $$T_{G_i^c} := \min\{n > 1 | X_n \notin G_i\}, \qquad i = 1, 2, \ldots, d.$$

Since the state $z$ in Assumption 2.1 is a minimizer of $J^*(\lambda, \cdot)$, by Theorem 4.4, it follows that $z \notin G_i$ when $1 \leq i \leq d$, by (4.5) and (5.1). Therefore, (5.3) and (2.4) together imply that

(5.4) $$T_{G_i^c} \leq T$$

and, via Lemma 4.2, this yields

(5.5) $$P_x^\pi[T_{G_i^c} = n] \leq P_x^\pi[T_{G_i^c} \geq n] \leq P_x^\pi[T \geq n] \leq \beta_0 \beta^n,$$
$$n \in \mathbb{N}, i = 1, 2, \ldots, d.$$

The proof of Theorem 5.1, which parallels the ideas used to establish Theorem 4.4(ii), relies on the following lemma extending conclusions in Lemmas 4.3(iii) and 4.8.

LEMMA 5.2. *Let $\varepsilon \in (0, \xi)$ and $x \in G_i$ be arbitrary but fixed, where $i > 0$, and suppose that $\pi \in \mathcal{P}$ is $\varepsilon$-optimal at state $x$. In this case, assertions* (i)–(iv) *below are valid.*

(i) $E_x^\pi[\exp\{\lambda \alpha \sum_{t=0}^{T_{G_i^c}-1} [C(X_t, A_t) - J^*(\lambda, X_t)]\}] < \infty.$
(ii) $P_x^\pi[J^*(\lambda, X_r) \leq \gamma_i] = 1$ *for each* $r \in \mathbb{N}.$

*Consequently:*

(iii) *When $X_0 = x$ and the system is driven by $\pi$, the inclusion $A_t \in B^*(X_t)$ holds before $T_{G_i^c}$ with probability 1, that is,*

$$P_x^\pi\left[\bigcap_{t=0}^{T_{G_i^c}-1} [A_t \in B^*(X_t)]\right] = 1.$$

(iv) $P_x^\pi[X_{T_{G_i^c}} \in \bigcup_{k=0}^{i-1} G_k] = 1.$

PROOF. (i) The argument is along the lines in the proof of Lemma 4.3(iii). First, notice that (5.3) yields that $X_t \in G_i$ if $1 \leq t < T_{G_i^c}$, and then, $J^*(\lambda, X_t) = \gamma_i$ for $0 \leq t < T_{G_i^c}$ when $X_0 \in G_i$. Therefore, using that $x \in G_i$, via Hölder's inequality it follows that

$$E_x^\pi\left[\exp\left\{\lambda \alpha \sum_{t=0}^{T_{G_i^c}-1} [C(X_t, A_t) - J^*(\lambda, X_t)]\right\}\right]$$
$$= E_x^\pi\left[\exp\left\{\lambda \alpha \sum_{t=0}^{T_{G_i^c}-1} [C(X_t, A_t) - \gamma_i]\right\}\right]$$

$$= \sum_{n=1}^{\infty} E_x^\pi \left[ \exp\left\{ \lambda\alpha \sum_{t=0}^{n-1} [C(X_t, A_t) - \gamma_i] \right\} I[T_{G_i^c} = n] \right]$$

$$\leq \sum_{n=1}^{\infty} \left( E_x^\pi \left[ \exp\left\{ \lambda \sum_{t=0}^{n-1} [C(X_t, A_t) - \gamma_i] \right\} \right] \right)^\alpha (P_x^\pi[T_{G_i^c} = n])^{1-\alpha}$$

so that

(5.6)
$$E_x^\pi \left[ \exp\left\{ \lambda\alpha \sum_{t=0}^{T_{G_i^c}-1} [C(X_t, A_t) - J^*(\lambda, X_t)] \right\} \right]$$
$$\leq \sum_{n=1}^{\infty} e^{\lambda\alpha[J_n(\lambda,\pi,x) - n\gamma_i]} (P_x^\pi[T_{G_i^c} = n])^{1-\alpha};$$

see (2.1). Since $J^*(\lambda, x) = \gamma_i$ and $\pi$ is $\varepsilon$-optimal at $x$, it follows that, for some positive integer $n_0$, $J_n(\lambda, \pi, x) \leq n(\gamma_i + \varepsilon)$ when $n > n_0$. This leads, via (5.5), to

$$e^{\lambda\alpha[J_n(\lambda,\pi,x) - n\gamma_i]} (P_x^\pi[T_{G_i^c} = n])^{1-\alpha} \leq \beta_0^{1-\alpha} (e^{\lambda\alpha\varepsilon} \beta^{1-\alpha})^n, \qquad n > n_0.$$

Observing that the inclusion $\varepsilon \in (0, \xi)$ yields that $e^{\lambda\alpha\varepsilon} \beta^{1-\alpha} < 1$ [see (4.2) and Definition 4.6], the above-displayed inequality and (5.6) together imply that $E_x^\pi[\exp\{\lambda\alpha \sum_{t=0}^{T_{G_i^c}-1} [C(X_t, A_t) - J^*(\lambda, X_t)]\}]$ is finite.

(ii) Let $r \in \mathbb{N}$ be arbitrary but fixed. Since policy $\pi$ is $\varepsilon$-optimal at $x$, Lemma 4.5(ii) yields that $P_x^\pi[J^*(\lambda, X_r) \leq J^*(\lambda, x) + \varepsilon] = 1$, and using the inclusion $\varepsilon \in (0, \xi)$, Remark 4.7 allows us to write $P_x^\pi[J^*(\lambda, X_r) \leq J^*(\lambda, x)] = 1$. The conclusion follows since $J^*(\lambda, x) = \gamma_i$.

(iii) Let $r, k \in \mathbb{N}$ be such that $r < k$, and suppose that the pair $(w, a) \in \mathbb{K}$ is such that

(5.7) $$P_x^\pi[X_r = w, A_r = a, T_{G_i^c} = k] > 0.$$

Since $r < k$, from the definition of the exit time $T_{G_i^c}$ and the inclusion $x \in G_i$, it follows that

$$w \in G_i,$$

and it will be shown that $a \in B^*(w)$. To achieve this goal, suppose that $y \in S$ satisfies $p_{wy}(a) > 0$ and observe that the Markov property yields $P_x^\pi[X_{r+1} = y | X_r = w, A_r = a] = p_{wy}(a) > 0$; since $P_x^\pi[X_r = w, A_r = a] > 0$, by (5.7), it follows that

$$P_x^\pi[X_{r+1} = y] \geq P_x^\pi[X_{r+1} = y, X_r = w, A_r = a]$$
$$\geq P_x^\pi[X_{r+1} = y | X_r = w, A_r = a] P_x^\pi[X_r = w, A_r = a] > 0.$$

Therefore, part (ii) yields that $J^*(\lambda, y) \leq \gamma_i$. Since $y \in S$ satisfying $p_{wy}(a) > 0$ is arbitrary, it follows that

$$\max\{J^*(\lambda, y)|p_{wy}(a) > 0\} \leq \gamma_i = J^*(\lambda, w),$$

where the inclusion $w \in G_i$ was used to set the equality. Then

$$\max\{J^*(\lambda, y)|p_{wy}(a) > 0\} = J^*(\lambda, w),$$

by Lemma 3.1, so that $a \in B^*(w)$; see Definition 4.1(i). In short, when $r < k$,

$$P_x^\pi[X_r = w, A_r = a, T_{G_i^c} = k] > 0 \implies a \in B^*(w),$$

and it follows that

$$P_x^\pi[T_{G_i^c} = k] = \sum_{(w,a) \in \mathbb{K}} P_x^\pi[X_r = w, A_r = a, T_{G_i^c} = k]$$

$$= \sum_{(w,a) \in \mathbb{K}, a \in B^*(w)} P_x^\pi[X_r = w, A_r = a, T_{G_i^c} = k]$$

and then

$$P_x^\pi[T_{G_i^c} = k] = P_x^\pi[[A_r \in B^*(X_r)] \cap [T_{G_i^c} = k]].$$

Since this equality holds whenever $r < k$, it follows that

$$P_x^\pi[T_{G_i^c} = k] = P_x^\pi\left[\bigcap_{r=0}^{k-1}[A_r \in B^*(X_r)] \cap [T_{G_i^c} = k]\right]$$

$$= P_x^\pi\left[\bigcap_{r=0}^{T_{G_i^c}-1}[A_r \in B^*(X_r)] \cap [T_{G_i^c} = k]\right].$$

Summing up over all positive integers $k$, this yields

$$P_x^\pi[T_{G_i^c} < \infty] = P_x^\pi\left[\bigcap_{r=0}^{T_{G_i^c}-1}[A_r \in B^*(X_r)] \cap [T_{G_i^c} < \infty]\right],$$

and the conclusion follows since, by (5.5), $P_x^\pi[T_{G_i^c} < \infty] = 1$.

(iv) Notice that (4.5), (5.1) and part (ii) together yield that, for each positive integer $r$,

$$P_x^\pi\left[X_r \in \bigcup_{k=0}^{i} G_k\right] = 1.$$

On the other hand, from (5.3) it follows that $X_r \notin G_i$ on the event $[T_{G_i^c} = r]$, so that the above displayed equation implies that

$$P_x^\pi\left[T_{G_i^c} = r, X_{T_{G_i^c}} \in \bigcup_{k=0}^{i-1} G_k\right] = P_x^\pi\left[T_{G_i^c} = r, X_r \in \bigcup_{k=0}^{i-1} G_k\right] = P_x^\pi[T_{G_i^c} = r].$$

Hence,

$$P_x^\pi\left[T_{G_i^c} < \infty, X_{T_{G_i^c}} \in \bigcup_{k=0}^{i-1} G_k\right] = \sum_{r=1}^{\infty} P_x^\pi\left[T_{G_i^c} = r, X_{T_{G_i^c}} \in \bigcup_{k=0}^{i-1} G_k\right]$$

$$= \sum_{r=1}^{\infty} P_x^\pi[T_{G_i^c} = r] = P_x^\pi[T_{G_i^c} < \infty],$$

and the conclusion follows using that $T_{G_i^c}$ is finite with probability 1. □

PROOF OF THEOREM 5.1. For each $m = 0, 1, 2, \ldots, d$, consider the following claim:

($C_m$) $h(x) < \infty$ for every $x \in G_m$.

Observe that the conclusion of Theorem 5.1 is equivalent to the truth of every ($C_m$) a fact that will be established by induction. To begin with, notice that ($C_0$) holds, by Theorem 4.4(ii); see (4.5) and (5.1). Assume now that $i \leq d$ is a positive integer such that ($C_m$) holds when $m < i$. Under this condition it will be proved that ($C_i$) is valid. From this induction hypothesis, the definition of $h(\cdot)$ in (4.1) yields that for each $y \in \bigcup_{k=0}^{i-1} G_k$ there exists a policy $\delta^y$ such that

$$(5.8) \quad \delta^y \in \mathcal{P}^* \quad \text{and} \quad E_y^{\delta^y}\left[\exp\left\{\lambda\alpha \sum_{t=0}^{T-1}[C(X_t, A_t) - J^*(\lambda, X_t)]\right\}\right] < \infty.$$

Next, let $x \in G_i$ and $\varepsilon \in (0, \xi)$ be arbitrary but fixed, and select a policy $\pi \in \mathcal{P}$ which is $\varepsilon$-optimal at $x$. Given a stationary policy $f$ satisfying $f(w) \in B^*(w)$ for every $w \in S$, define the new policy $\delta$ as follows: For each $t \in \mathbb{N}$ and $\mathbf{h}_t \in \mathbb{H}_t$:

(a) If $x_0 \neq x$, $\delta_t(\{f(x_t)\}|\mathbf{h}_t) = 1$.
(b) If $x_0 = x$ and $x_k \in G_i$ for every $k = 1, 2, \ldots, t$, then $\delta_t(\cdot|\mathbf{h}_t) = \pi_t(\cdot|\mathbf{h}_t)$.
(c) If $x_0 = x$ and, for some positive integer $r \leq t$, $x_k \in G_i$ for every $k < r$ and $x_r \notin G_i$, then

$$\delta_t(\cdot|\mathbf{h}_t) = \delta_{t-r}^{x_r}(\cdot|x_r, \ldots, x_{t-1}, a_{t-1}, x_t) \quad \text{if } x_r \in \bigcup_{j=0}^{i-1} G_j,$$

$$\delta_t(\{f(x_t)\}|\mathbf{h}_t) = 1 \quad \text{when } x_r \in S \setminus \bigcup_{j=0}^{i} G_j.$$

A controller driving the system via policy $\delta$ operates as follows: When the initial state is $X_0 \neq x$, at each decision time the actions are selected according to $f$. On the other hand, when $X_0 = x$, she uses policy $\pi$ to choose actions while the system stays in $G_i$, but when the system first leaves $G_i$ at time $T_{G_i^c} = k$, then the decision maker "forgets" the history observed before time $k$ and, as if the process had started again, she switches to policy $\delta^{X_k}$ if $X_k$ belongs to some set $G_m$ with $m < i$, or to policy $f$ otherwise. Now, let $t \in \mathbb{N}$ be fixed. When the initial state is $x$, $\delta$ and $\pi$ coincide while the system stays in $G_i$, by part (b) so that Lemma 5.2(iii) yields that

$$A_t \in B^*(X_t)$$

holds on $[T_{G_i^c} < t]$ $P_x^\delta$-a.s. whereas, by part (c), the choice of $f$ and the inclusion in (5.8) imply that the above displayed relation also occurs $P_x^\delta$-a.s. on the event $[T_{G_i^c} \geq t]$. When $X_0 = w \neq x$, from the choice of $f$ and part (a) in the above definition, it follows that $P_w^\delta[A_t \in B^*(X_t)] = 1$, so that

(5.9) $$\delta \in \mathcal{P}^*;$$

see Definition 4.1. Moreover, using again that $\delta$ and $\pi$ coincide before $T_{G_i^c}$ when $x$ is the initial state, it follows that the event $[X_{T_{G_i^c}} \in \bigcup_{k=0}^{i-1} G_k]$ has the same probability with respect to $P_x^\delta$ and $P_x^\pi$, whereas the expectation of $\exp\{\lambda\alpha \sum_{t=0}^{T_{G_i^c}-1}[C(X_t, A_t) - \gamma_i]\}$ with respect to these measures coincides. Thus, by parts (i) and (iv) of Lemma 5.2,

(5.10) $$P_x^\delta\left[X_{T_{G_i^c}} \in \bigcup_{k=0}^{i-1} G_k\right] = 1$$

and

(5.11) $$E_x^\delta\left[\exp\left\{\lambda\alpha \sum_{t=0}^{T_{G_i^c}-1}[C(X_t, A_t) - J^*(\lambda, X_t)]\right\}\right] < \infty.$$

Next, it will be shown that

(5.12) $$E_x^\delta\left[\exp\left\{\lambda\alpha \sum_{t=0}^{T-1}[C(X_t, A_t) - J^*(\lambda, X_t)]\right\}\right] < \infty.$$

To achieve this goal, notice that

$$E_x^\delta\left[\exp\left\{\lambda\alpha \sum_{t=0}^{T-1}[C(X_t, A_t) - J^*(\lambda, X_t)]\right\}I[T_{G_i^c} = T]\right]$$

$$= E_x^\delta\left[\exp\left\{\lambda\alpha \sum_{t=0}^{T_{G_i^c}-1}[C(X_t, A_t) - J^*(\lambda, X_t)]\right\}I[T_{G_i^c} = T]\right]$$

$$\leq E_x^\delta\left[\exp\left\{\lambda\alpha \sum_{t=0}^{T_{G_i^c}-1}[C(X_t, A_t) - J^*(\lambda, X_t)]\right\}\right]$$

and then

$$(5.13) \quad E_x^\delta\left[\exp\left\{\lambda\alpha \sum_{t=0}^{T-1}[C(X_t, A_t) - J^*(\lambda, X_t)]\right\}I[T_{G_i^c} = T]\right] < \infty,$$

by (5.11). Next, observe that for each positive integer $r$,

$$E_x^\delta\left[\exp\left\{\lambda\alpha \sum_{t=0}^{T-1}[C(X_t, A_t) - J^*(\lambda, X_t)]\right\}I[r = T_{G_i^c} < T]|I_r\right]$$

$$= \exp\left\{\lambda\alpha \sum_{t=0}^{r-1}[C(X_t, A_t) - J^*(\lambda, X_t)]\right\}I[r = T_{G_i^c} < T]$$

$$\times E_x^\delta\left[\exp\left\{\lambda\alpha \sum_{t=r}^{T-1}[C(X_t, A_t) - J^*(\lambda, X_t)]\right\}|I_r\right]$$

$$= \exp\left\{\lambda\alpha \sum_{t=0}^{T_{G_i^c}-1}[C(X_t, A_t) - J^*(\lambda, X_t)]\right\}I[r = T_{G_i^c} < T]$$

$$\times E_x^\delta\left[\exp\left\{\lambda\alpha \sum_{t=r}^{T-1}[C(X_t, A_t) - J^*(\lambda, X_t)]\right\}|I_r\right]$$

and that, on the event $[T_{G_i^c} = r]$, $X_r$ lies in $\bigcup_{j=0}^{r-1} G_i$ $P_x^\delta$-a.s., by (5.10). Thus, part (c) in the definition of policy $\delta$ yields, via the Markov property, that the following holds with probability 1 with respect to $P_x^\delta$:

$$I[r = T_{G_i^c} < T]E_x^\delta\left[\exp\left\{\lambda\alpha \sum_{t=r}^{T-1}[C(X_t, A_t) - J^*(\lambda, X_t)]\right\}|I_r\right]$$

$$= I\left[r = T_{G_i^c} < T, X_r \in \bigcup_{j=0}^{i-1} G_i\right]$$

$$\times E_{X_r}^{\delta^{X_r}}\left[\exp\left\{\lambda\alpha \sum_{t=0}^{T-1}[C(X_t, A_t) - J^*(\lambda, X_t)]\right\}\right]$$

$$\leq MI[r = T_{G_i^c} < T],$$

where $M := \max\{E_y^{\delta^y}[\exp\{\lambda\alpha \sum_{t=0}^{T-1}[C(X_t) - V(X_t)]\}]|y \in \bigcup_{j=0}^{i-1} G_i\} < \infty$, and the inequality is due to the induction hypothesis. Combining the last two displayed relations, it follows that

$$E_x^\delta\left[\exp\left\{\lambda\alpha \sum_{t=0}^{T-1}[C(X_t, A_t) - J^*(\lambda, X_t)]\right\}I[r = T_{G_i^c} < T]|I_r\right]$$

$$\leq MI[r = T_{G_i^c} < T] \exp\left\{\lambda\alpha \sum_{t=0}^{T_{G_i^c}-1} [C(X_t, A_t) - J^*(\lambda, X_t)]\right\}, \qquad P_x^\delta\text{-a.s.,}$$

so that

$$E_x^\delta\left[\exp\left\{\lambda\alpha \sum_{t=0}^{T-1}[C(X_t, A_t) - J^*(\lambda, X_t)]\right\}I[r = T_{G_i^c} < T]\right]$$

$$\leq ME_x^\delta\left[I[r = T_{G_i^c} < T] \exp\left\{\lambda\alpha \sum_{t=0}^{T_{G_i^c}-1} [C(X_t, A_t) - J^*(\lambda, X_t)]\right\}\right].$$

Since this inequality is valid for every positive integer $r$ and $T_{G_i^c}$ is finite $P_x^\delta$-a.s., it follows that

$$E_x^\delta\left[\exp\left\{\lambda\alpha \sum_{t=0}^{T-1}[C(X_t, A_t) - J^*(\lambda, X_t)]\right\}I[T_{G_i^c} < T]\right]$$

$$\leq ME_x^\delta\left[I[T_{G_i^c} < T] \exp\left\{\lambda\alpha \sum_{t=0}^{T_{G_i^c}-1} [C(X_t, A_t) - J^*(\lambda, X_t)]\right\}\right]$$

$$\leq ME_x^\delta\left[\exp\left\{\lambda\alpha \sum_{t=0}^{T_{G_i^c}-1} [C(X_t, A_t) - J^*(\lambda, X_t)]\right\}\right].$$

Since $M$ is finite, this relation and (5.11) yield that $E_x^\delta[\exp\{\lambda\alpha \sum_{t=0}^{T-1}[C(X_t, A_t) - J^*(\lambda, X_t)]\}I[T_{G_i^c} < T]]$ is finite, and this fact, (5.4) and (5.13) together imply that (5.12) occurs. To conclude, observe that, by the definition of function $h(\cdot)$ in (4.1), the inclusion in (5.9) and (5.12) together yield that $h(x) < \infty$; since $h(\cdot) > -\infty$, by Lemma 4.3, it follows that $h(x)$ is finite and, since $x \in G_i$ is arbitrary, this shows that claim $(C_i)$ holds, completing the induction proof. □

**6. A key inequality.** This section contains the last technical tool that, together with the finiteness result in Theorem 5.1, will be used to establish Theorem 3.5. The main objective is to establish the following.

THEOREM 6.1. *Let $z$ be the state in Assumption 2.1. In this case, the deviation function in (4.1) satisfies that*

$$h(z) \leq 0.$$

The proof of this theorem relies on Lemma 6.3, whose conclusions involve the random times at which the system occupies the distinguished state $z$ in Assumption 2.1.

DEFINITION 6.2. (i) The sequence $\{T_k\}$ of successive arrival times to state $z$ is recursively determined as follows:

$$T_1 := T \quad \text{and} \quad T_k := \min\{n > T_{k-1} | X_n = z\}, \qquad k > 1;$$

see (2.4) for the definition of $T$.

(ii) Given $\varepsilon > 0$, define $\psi(\varepsilon)$ by

$$\psi(\varepsilon) := \frac{1}{\lambda} \inf_{\delta \in \mathcal{P}^*} \log \left( E_z^\delta \left[ \exp\left\{ \lambda \sum_{t=0}^{T-1} [C(X_t, A_t) - \gamma_0 - 2\varepsilon] \right\} \right] \right)$$

where, as before, $\gamma_0 = J^*(\lambda, z)$.

It is not difficult to see that each $T_k$ is a stopping time with respect to the family of $\sigma$-fields $\{\sigma(I_n)\}$, that is, the event $[T_k = m]$ always lies in $\sigma(I_m)$, and that

(6.1) $$T_k \geq k, \qquad k = 1, 2, \ldots.$$

LEMMA 6.3. *Let $\varepsilon \in (0, \xi)$ be fixed, and suppose that $\pi \in \mathcal{P}$ is $\varepsilon$-optimal at state $z$. In this case:*

(i) *For each positive integer $k$,*

(6.2) $$e^{\lambda k \psi(\varepsilon)} \leq E_z^\pi \left[ \exp\left\{ \lambda \sum_{t=0}^{T_k-1} [C(X_t, A_t) - \gamma_0 - 2\varepsilon] \right\} \right].$$

(ii) *There exists $n_0$ such that, for every $k \geq n_0$,*

$$E_z^\pi \left[ \exp\left\{ \lambda \sum_{t=0}^{T_k-1} [C(X_t, A_t) - \gamma_0 - 2\varepsilon] \right\} \right] \leq \frac{e^{-\lambda \varepsilon k}}{1 - e^{-\lambda \varepsilon}}.$$

*Consequently:*

(iii) $\psi(\varepsilon) \leq -\varepsilon.$

PROOF. (i) The argument is by induction. Since $\pi$ is $\varepsilon$-optimal at $z$ and $\varepsilon \in (0, \xi)$, Theorem 4.4(i) and Lemma 4.8(iii) together yield that there exists $\delta \in \mathcal{P}^*$ such that $P_z^\pi = P_z^\delta$. In this case, using that $T_1 = T$, it follows that

$$E_z^\pi \left[ \exp\left\{ \lambda \sum_{t=0}^{T_1-1} [C(X_t, A_t) - \gamma_0 - 2\varepsilon] \right\} \right]$$

$$= E_z^\delta \left[ \exp\left\{ \lambda \sum_{t=0}^{T-1} [C(X_t, A_t) - \gamma_0 - 2\varepsilon] \right\} \right] \geq e^{\lambda \psi(\varepsilon)},$$

where the inequality is due to Definition 6.2(ii), so that (6.2) is valid for $k = 1$. Let $n$ be an integer larger than 1, and suppose that (6.2) holds when $k < n$.

In this situation, take a positive integer $r$ and $\tilde{\mathbf{h}}_r = (\tilde{x}_0, \tilde{a}_0, \ldots, \tilde{x}_{r-1}, \tilde{a}_{r-1}, \tilde{x}_r) \in \mathbb{H}_r$ satisfying that

(6.3) $$P_x^\pi[T_{n-1} = r, I_r = \tilde{\mathbf{h}}_r] > 0;$$

since $X_{T_{n-1}} = z$ when $T_{n-1}$ is finite, it follows that $\tilde{x}_r = z$. Notice now that $T_n > T_{n-1}$ so that, on the event $[T_{n-1} = r]$,

$$\sum_{t=0}^{T_n-1} [C(X_t, A_t) - \gamma_0 - 2\varepsilon]$$
$$= \sum_{t=0}^{r-1} [C(X_t, A_t) - \gamma_0 - 2\varepsilon] + \sum_{t=r}^{T_n-1} [C(X_t, A_t) - \gamma_0 - 2\varepsilon]$$

and an application of the Markov property yields, via Definition 6.2, that

$$E_z^\pi\left[\exp\left\{\lambda \sum_{t=0}^{T_n-1}[C(X_t, A_t) - \gamma_0 - 2\varepsilon]\right\}\bigg| T_{n-1} = r, I_r = \tilde{\mathbf{h}}_r\right]$$
$$= \exp\left\{\lambda \sum_{t=0}^{r-1}[C(X_t, A_t) - \gamma_0 - 2\varepsilon]\right\}$$
$$\times E_z^\delta\left[\exp\left\{\lambda \sum_{t=0}^{T_1}[C(X_t, A_t) - \gamma_0 - 2\varepsilon]\right\}\right],$$

where the shifted policy $\delta$ is as in (4.3). Since $\pi$ is $\varepsilon$-optimal at $z$, Lemma 4.5(i) yields that $J(\lambda, \delta, z) \leq J(\lambda, \pi, z) \leq J^*(\lambda, z) + \varepsilon$, so that $\delta$ itself is $\varepsilon$-optimal at $z$. Therefore, applying the case $k = 1$ of (6.2) to this policy $\delta$, it follows that $E_z^\delta[\exp\{\lambda \sum_{t=0}^{T_1}[C(X_t, A_t) - \gamma_0 - 2\varepsilon]\}] \geq e^{\lambda \psi(\varepsilon)}$, which combined with the above displayed equation leads to

$$E_z^\pi\left[\exp\left\{\lambda \sum_{t=0}^{T_n-1}[C(X_t, A_t) - \gamma_0 - 2\varepsilon]\right\}\bigg| T_{n-1} = r, I_r = \tilde{\mathbf{h}}_r\right]$$
$$\geq \exp\left\{\lambda \sum_{t=0}^{r-1}[C(X_t, A_t) - \gamma_0 - 2\varepsilon]\right\} e^{\lambda \psi(\varepsilon)}.$$

Since this inequality is valid whenever (6.3) holds, it follows that

$$E_z^\pi\left[\exp\left\{\lambda \sum_{t=0}^{T_n-1}[C(X_t, A_t) - \gamma_0 - 2\varepsilon]\right\}\right]$$
$$\geq E_z^\pi\left[\exp\left\{\lambda \sum_{t=0}^{T_{n-1}-1}[C(X_t, A_t) - \gamma_0 - 2\varepsilon]\right\}\right] e^{\lambda \psi(\varepsilon)} \geq e^{n\lambda \psi(\varepsilon)},$$

where the induction hypothesis was used to set the second inequality. This establishes the case $k = n$ of (6.2) and completes the induction argument.

(ii) Since $\pi$ is $\varepsilon$-optimal at $z$, there exists a positive integer $n_0$ such that $J_n(\lambda, \pi, z) \leq n(J^*(\lambda, z) + \varepsilon) = n(\gamma_0 + \varepsilon)$ when $n \geq n_0$. Observing that

$$E_z^\pi \left[ \exp\left\{ \lambda \sum_{t=0}^{n-1} [C(X_t, A_t) - \gamma_0 - 2\varepsilon] \right\} \right]$$
$$= e^{-2n\lambda\varepsilon - n\gamma_0} E_z^\pi \left[ \exp\left\{ \lambda \sum_{t=0}^{n-1} C(X_t, A_t) \right\} \right]$$
$$= e^{-2n\lambda\varepsilon - n\gamma_0} e^{J_n(\lambda, \pi, z)},$$

it follows that

(6.4) $\quad E_z^\pi \left[ \exp\left\{ \lambda \sum_{t=0}^{n-1} [C(X_t, A_t) - \gamma_0 - 2\varepsilon] \right\} \right] \leq e^{-n\lambda\varepsilon}, \qquad n \geq n_0.$

Next, let the positive integer $k$ and $\rho \in (0,1)$ be fixed. In this case, (6.1) and Hölder's inequality yield that

$$E_z^\pi \left[ \exp\left\{ \lambda\rho \sum_{t=0}^{T_k-1} [C(X_t, A_t) - \gamma_0 - 2\varepsilon] \right\} \right]$$
$$= \sum_{n=k}^\infty E_z^\pi \left[ \exp\left\{ \lambda\rho \sum_{t=0}^{n-1} [C(X_t, A_t) - \gamma_0 - 2\varepsilon] \right\} I[T_k = n] \right]$$
$$\leq \sum_{n=k}^\infty \left( E_z^\pi \left[ \exp\left\{ \lambda \sum_{t=0}^{n-1} [C(X_t, A_t) - \gamma_0 - 2\varepsilon] \right\} \right] \right)^\rho (P_x^\pi I[T_k = n])^{1-\rho}$$
$$\leq \sum_{n=k}^\infty \left( E_z^\pi \left[ \exp\left\{ \lambda \sum_{t=0}^{n-1} [C(X_t, A_t) - \gamma_0 - 2\varepsilon] \right\} \right] \right)^\rho.$$

Combining this with (6.4) it follows that

$$E_z^\pi \left[ \exp\left\{ \lambda\rho \sum_{t=0}^{T_k-1} [C(X_t, A_t) - \gamma_0 - 2\varepsilon] \right\} \right] \leq \sum_{n=k}^\infty e^{-n\varepsilon\rho\lambda} = \frac{e^{-k\varepsilon\lambda\rho}}{1 - e^{-\varepsilon\lambda\rho}},$$
$$k \geq n_0.$$

Given a sequence $\{\rho_m\}$ of positive numbers increasing to 1, this inequality implies, via Fatou's lemma, that for every positive integer $k \geq n_0$,

$$E_z^\pi \left[ \exp\left\{ \lambda \sum_{t=0}^{T_k-1} [C(X_t, A_t) - \gamma_0 - 2\varepsilon] \right\} \right]$$
$$= E_z^\pi \left[ \liminf_{m \to \infty} \exp\left\{ \lambda\rho_m \sum_{t=0}^{T_k-1} [C(X_t, A_t) - \gamma_0 - 2\varepsilon] \right\} \right]$$

$$\leq \liminf_{m \to \infty} E_z^\pi \left[ \exp\left\{ \lambda \rho_m \sum_{t=0}^{T_k-1} [C(X_t, A_t) - \gamma_0 - 2\varepsilon] \right\} \right]$$

$$\leq \liminf_{m \to \infty} \frac{e^{-k\varepsilon \lambda \rho_m}}{1 - e^{-\varepsilon \lambda \rho_m}}$$

and then

$$E_z^\pi \left[ \exp\left\{ \lambda \sum_{t=0}^{T_k-1} [C(X_t, A_t) - \gamma_0 - 2\varepsilon] \right\} \right] \leq \frac{e^{-k\varepsilon \lambda}}{1 - e^{-\varepsilon \lambda}}, \qquad k \geq n_0.$$

(iii) Observe that parts (i) and (ii) together yield that $e^{\lambda k \psi(\varepsilon)} \leq e^{-k\varepsilon \lambda}/(1 - e^{-\varepsilon \lambda})$ when $k$ is large enough, and in this case

$$\psi(\varepsilon) \leq -\varepsilon - \frac{1}{\lambda k} \log(1 - e^{\varepsilon \lambda}),$$

so that the conclusion follows letting $k$ increase to $\infty$. $\square$

PROOF OF THEOREM 6.1. It will be shown that there exists a policy $\delta$ satisfying

(6.5) $\qquad \delta \in \mathcal{P}^* \quad \text{and} \quad E_z^\delta \left[ \exp\left\{ \lambda \sum_{t=0}^{T-1} [C(X_t, A_t) - \gamma_0] \right\} \right] \leq 1.$

Assuming that such a policy exists, Theorem 6.1 can be established as follows: First, recall that $J^*(\lambda, \cdot)$ satisfies the min–max equation in (3.3), by Lemma 3.1, and that $B^* = B_{J^*(\lambda, \cdot)}$, by Definition 4.1(i). Therefore, the inclusion $\delta \in \mathcal{P}^*$ yields that $P_z^\delta[A_t \in B_{J^*(\lambda, \cdot)}(X_t)] = 1$ for every $t \in \mathbb{N}$, by Definition 4.1(ii), so that an application of Lemma 3.4(i) implies that, for each $n \in \mathbb{N}$,

$$J^*(\lambda, X_n) \leq J^*(\lambda, X_0) = J^*(\lambda, z) = \gamma_0, \qquad P_z^\delta\text{-a.s.}$$

Since $\gamma_0$ is the minimum value of $J^*(\lambda, \cdot)$, it follows that $P_z^\delta[J^*(\lambda, X_n) = \gamma_0] = 1$ for every $n \in \mathbb{N}$, so that the inequality in (6.5) is equivalent to

$$E_z^\delta \left[ \exp\left\{ \lambda \sum_{t=0}^{T-1} [C(X_t, A_t) - J^*(\lambda, X_t)] \right\} \right] \leq 1.$$

From this point, an application of Hölder's inequality yields that

$$E_z^\delta \left[ \exp\left\{ \lambda \alpha \sum_{t=0}^{T-1} [C(X_t, A_t) - J^*(\lambda, X_t)] \right\} \right]$$

$$\leq \left( E_z^\delta \left[ \exp\left\{ \lambda \sum_{t=0}^{T-1} [C(X_t, A_t) - J^*(\lambda, X_t)] \right\} \right] \right)^\alpha \leq 1;$$

recall that the fixed number $\alpha$ lies in $(0,1)$. Combining this inequality with the inclusion $\delta \in \mathcal{P}^*$ and (4.1), it follows that

$$h(z) \leq \frac{1}{\lambda} \log\left( E_z^\delta\left[ \exp\left\{ \lambda\alpha \sum_{t=0}^{T-1} [C(X_t, A_t) - J^*(\lambda, X_t)] \right\} \right] \right) \leq 0,$$

completing the proof of Theorem 6.1. To conclude, (6.5) will be established. Let $\{\varepsilon_k\} \subset (0, \xi)$ be a sequence converging to zero and notice that, for each $k \in \mathbb{N}$, Lemma 6.3(iii) yields that $e^{\lambda\psi(\varepsilon_k)} \leq e^{-\lambda\varepsilon_k} < e^{-\lambda\varepsilon_k/2}$. Thus, by Definition 6.2, for every $k \in \mathbb{N}$ there exists a policy $\pi^k \in \mathcal{P}^*$ such that

$$(6.6) \qquad E_z^{\pi^k}\left[ \exp\left\{ \lambda \sum_{t=0}^{T-1} [C(X_t, A_t) - \gamma_0 - 2\varepsilon_k] \right\} \right] \leq e^{-\lambda\varepsilon_k/2}.$$

Let $\mathbb{P}(A)$ be the class of probability measures defined on the subsets of the action space $A$. For each $r \in \mathbb{N}$ and $\mathbf{h}_r \in \mathbb{H}_r$, $\{\pi_r^k(\cdot|\mathbf{h}_r)|k \in \mathbb{N}\}$ is a sequence in $\mathbb{P}(A)$ and, since $A$ is finite, there exists $\delta_r(\cdot|\mathbf{h}_r) \in \mathbb{P}(A)$, as well as a subsequence of $\{\pi^k\}$, denoted by $\{\pi^m\}$, such that

$$(6.7) \qquad \lim_{m\to\infty} \pi_r^m(F|\mathbf{h}_r) =: \delta_r(F|\mathbf{h}_r), \qquad F \subset A.$$

Moreover, since $\bigcup_{r=0}^{\infty} \mathbb{H}_r$ is denumerable, applying Cantor's diagonal method it can be assumed that this convergence holds for every $F \subset A$, $r \in \mathbb{N}$ and $\mathbf{h}_r \in \mathbb{H}_r$, and it will be shown that $\delta := \{\delta_r\}$ satisfies (6.5). To achieve this goal, first notice that $\pi_r^k(A(x_r)|\mathbf{h}_r) = 1$ always holds, since $\pi^k \in \mathcal{P}^* \subset \mathcal{P}$, so that (6.7) yields that $\delta_r(A(x_r)|\mathbf{h}_r) = 1$ for every $r \in \mathbb{N}$ and $\mathbf{h}_r \in \mathbb{H}_r$, that is, $\delta$ is a policy. Next, observe that the equality

$$P_x^\pi[I_r = \mathbf{h}_r] = \delta_{x,x_0}\pi_0(a_0|x_0)p_{x_0x_1}(a_0)\pi_1(a_1|x_0, a_0, x_1) \times \cdots$$
$$\times \pi_{r-1}(a_{r-1}|x_0, a_0, \ldots, x_{n-1})p_{x_{n-1}x_n}(a_{n-a})$$

is always valid, where $\delta_{x,y} := 1$ if $x = y$ and $\delta_{x,y} := 0$ otherwise. Combining this equation with (6.7), it follows that for every $x \in S$, $r \in \mathbb{N}$ and $D: \mathbb{H}_r \to \mathbb{R}$,

$$(6.8) \qquad \lim_{m\to\infty} E_x^{\pi^m}[D(I_r)] = E_x^\delta[D(I_r)].$$

In particular, for each $n \in \mathbb{N}$ and $x \in S$, $P_x^\delta[A_n \in B^*(X_n)] = \lim_{m\to\infty} P_x^{\pi^m}[A_n \in B^*(X_n)] = 1$, where the inclusion $\pi^m \in \mathcal{P}^*$ was used to set the second equality, so that

$$(6.9) \qquad \delta \in \mathcal{P}^*;$$

see Definition 4.1. Moreover, (6.8) yields that for every $r \in \mathbb{N}$,

$$\lim_{m\to\infty} E_z^{\pi^m}\left[ \exp\left\{ \lambda \sum_{t=0}^{T-1} [C(X_t, A_t) - \gamma_0] \right\} I[T \leq r] \right]$$

(6.10)
$$= E_z^\delta\left[\exp\left\{\lambda \sum_{t=0}^{T-1}[C(X_t, A_t) - \gamma_0]\right\}I[T \leq r]\right].$$

Observing that
$$\exp\left\{\lambda \sum_{t=0}^{T-1}[C(X_t, A_t) - \gamma_0 - 2\varepsilon_k]\right\}I[T \leq r]$$
$$\geq e^{-2\lambda r \varepsilon_k} \exp\left\{\lambda \sum_{t=0}^{T-1}[C(X_t, A_t) - \gamma_0]\right\}I[T \leq r],$$

it follows that
$$E_z^{\pi^m}\left[\exp\left\{\lambda \sum_{t=0}^{T-1}[C(X_t, A_t) - \gamma_0]\right\}I[T \leq r]\right]$$
$$\leq e^{2\lambda r \varepsilon_m} E_z^{\pi^m}\left[\exp\left\{\lambda \sum_{t=0}^{T-1}[C(X_t, A_t) - \gamma_0 - 2\varepsilon_m]\right\}I[T \leq r]\right]$$
$$\leq e^{2\lambda r \varepsilon_m} E_z^{\pi^m}\left[\exp\left\{\lambda \sum_{t=0}^{T-1}[C(X_t, A_t) - \gamma_0 - 2\varepsilon_m]\right\}\right]$$

so that $E_z^{\pi^m}[\exp\{\lambda \sum_{t=0}^{T-1}[C(X_t, A_t) - \gamma_0]\}I[T \leq r]] \leq e^{2\lambda r \varepsilon_m + \lambda \varepsilon_m/2}$ [see (6.6)]; since $\{\varepsilon_m\}$ converges to zero, this inequality and (6.10) together yield that $E_z^\delta[\exp\{\lambda \sum_{t=0}^{T-1}[C(X_t, A_t) - \gamma_0]\}I[T \leq r]] \leq 1$ for every $r \in \mathbb{N}$ and, via the monotone convergence theorem, this implies that

$$E_z^\delta\left[\exp\left\{\lambda \sum_{t=0}^{T-1}[C(X_t, A_t) - \gamma_0]\right\}\right] \leq 1.$$

Combining this inequality with the inclusion in (6.9), it follows that the conditions in (6.5) are satisfied by policy $\delta$. □

**7. Proof of the main result.** After the previous preliminaries, in this section the characterization result in Theorem 3.5 will be finally proved. The argument combines Theorems 5.1 and 6.1 with the properties of the policies in $\mathcal{P}^*$ established in the following lemma.

LEMMA 7.1. *Given a policy $\pi \in \mathcal{P}^*$, suppose that for some $(x, a) \in \mathbb{K}$ the inequality $P_x^\pi[A_0 = a] > 0$ holds and define the shifted policy $\delta$ by*

(7.1) $$\delta_t(\cdot|\mathbf{h}_t) = \pi_{t+1}(\cdot|x, a, \mathbf{h}_t).$$

*In this case, for each $y \in S$ satisfying that $p_{xy}(a) > 0$, assertions* (i) *and* (ii) *below hold.*

(i) $P_y^\delta[A_t \in B^*(X_t)] = 1$ *for every* $t \in \mathbb{N}$.

(ii) *There exists* $\tilde\delta \in \mathcal{P}^*$ *such that* $P_y^\delta = P_y^{\tilde\delta}$.

PROOF. (i) Suppose that $p_{xy}(a) > 0$ and observe that $P_x^\pi[X_1 = y, A_0 = a] = P_x^\pi[X_1 = y | A_0 = a]P_x^\pi[A_0 = a] = p_{xy}(a)P_x^\pi[A_0 = a] > 0$. Thus, for every $t \in \mathbb{N}$,

$$P_x^\pi[A_{t+1} \notin B^*(X_{t+1})]$$
$$\geq P_x^\pi[A_{t+1} \notin B^*(X_{t+1}), X_1 = y, A_0 = a]$$
$$= P_x^\pi[A_{t+1} \notin B^*(X_{t+1}) | X_1 = y, A_0 = a]P_x^\pi[X_1 = y, A_0 = a].$$

Observing that $P_x^\pi[A_{t+1} \notin B^*(X_{t+1}) | X_1 = y, A_0 = a] = P_y^\delta[A_t \notin B^*(X_t)]$, which is due to the definition of policy $\delta$ and the Markov property, it follows that

$$P_x^\pi[A_{t+1} \notin B^*(X_{t+1})] \geq P_y^\delta[A_t \notin B^*(X_t)]P_x^\pi[X_1 = y, A_0 = a].$$

Since $P_x^\pi[A_{t+1} \notin B^*(X_{t+1})] = 0$, by the inclusion $\pi \in \mathcal{P}^*$, and $P_x^\pi[X_1 = y, A_0 = a] > 0$, it follows that $P_y^\delta[A_t \notin B^*(X_t)] = 0$, that is, $P_y^\delta[A_t \in B^*(X_t)] = 1$.

(ii) Pick a stationary policy $f$ such that $f(y) \in B^*(y)$ for each $y \in S$, and define the policy $\tilde\delta$ as follows: For each $t \in \mathbb{N}$ and $\mathbf{h}_t \in \mathbb{H}_t$,

$$\tilde\delta_t(\cdot|\mathbf{h}_t) = \delta_t(\cdot|\mathbf{h}_t) \quad \text{if } p_{xx_0}(a) > 0,$$
$$\tilde\delta_t(\{f(x_t)\}|\mathbf{h}_t) = 1 \quad \text{if } p_{xx_0}(a) = 0.$$

From this definition it follows that $P_w^{\tilde\delta} = P_w^\delta$ when $p_{xw}(a) > 0$, and $P_w^{\tilde\delta} = P_w^f$ if $p_{xw}(a) = 0$. Therefore, $P_y^{\tilde\delta} = P_y^\delta$, since $p_{xy}(a) > 0$, whereas the choice of $f$ and part (i) together imply that $P_y^{\tilde\delta}[A_t \in B^*(X_t)] = 1$ always holds, that is, $\tilde\delta \in \mathcal{P}^*$, by Definition 4.1. □

PROOF OF THEOREM 3.5. Recall that the fixed number $\alpha$ belongs to $(0, 1)$ and let $g(\cdot)$ be the function defined in (3.6). It will be shown that this function belongs to the family $\mathcal{G}$ in Definition 3.2. Using that $\alpha$ is positive, from Lemma 3.1 it is not difficult to see the min–max equation (3.3) holds, so that $g(\cdot)$ satisfies the first requirement in Definition 3.2. Moreover, for each $(x, a) \in \mathbb{K}$, the equality $g(x) = \max\{g(y)|p_{xy}(a) > 0\}$ is equivalent to $J^*(\lambda, x) = \max\{J^*(\lambda, y)|p_{xy}(a) > 0\}$, so that

$$B_g(x) = B_{J^*(\lambda,\cdot)} = B^*(x), \qquad x \in S;$$

see (3.5) and Definition 4.1(i). It will be verified that the second part of Definition 3.2 is satisfied by the pair $(g(\cdot), h(\cdot))$, where $h(\cdot)$ is given in (4.1). To achieve this goal, first notice that this function $h(\cdot)$ is finite, by Theorem 5.1. Next, select a policy $\pi \in \mathcal{P}^*$ and let $x \in S$ be arbitrary. For each action $a$ satisfying that $P_x^\pi[A_0 = a] > 0$, it follows that $a \in B^*(x)$, since $\pi \in \mathcal{P}^*$, whereas the Markov property yields

$$E_x^\pi\left[\exp\left\{\lambda\alpha \sum_{t=0}^{T-1}[C(X_t, A_t) - J^*(\lambda, X_t)]\right\}\Big| A_0 = a\right]$$

$$= E_x^\pi\left[\exp\left\{\lambda\alpha \sum_{t=0}^{T-1}[C(X_t, A_t) - J^*(\lambda, X_t)]\right\}I[T=1]\Big| A_0 = a\right]$$

$$+ E_x^\pi\left[\exp\left\{\lambda\alpha \sum_{t=0}^{T-1}[C(X_t, A_t) - J^*(\lambda, X_t)]\right\}I[T>1]\Big| A_0 = a\right]$$

$$= e^{\lambda\alpha[C(x,a) - J^*(\lambda,x)]}p_{xz}(a)$$

$$+ e^{\lambda\alpha[C(x,a) - J^*(\lambda,x)]}$$

$$\times \sum_{y\neq z} p_{xy}(a) E_y^\delta\left[\exp\left\{\lambda\alpha \sum_{t=0}^{T-1}[C(X_t, A_t) - J^*(\lambda, X_t)]\right\}\right],$$

where $\delta$ is the shifted policy in (7.1). By Lemma 7.1, there exists $\tilde{\delta} \in \mathcal{P}^*$ such that $P_y^\delta = P_y^{\tilde{\delta}}$ when $p_{xy}(a) > 0$, so that

$$E_x^\pi\left[\exp\left\{\lambda\alpha \sum_{t=0}^{T-1}[C(X_t, A_t) - J^*(\lambda, X_t)]\right\}\Big| A_0 = a\right]$$

$$= e^{\lambda\alpha[C(x,a) - J^*(\lambda,x)]}p_{xz}(a)$$

$$+ e^{\lambda\alpha[C(x,a) - J^*(\lambda,x)]}$$

$$\times \sum_{y\neq z} p_{xy}(a) E_y^{\tilde{\delta}}\left[\exp\left\{\lambda\alpha \sum_{t=0}^{T-1}[C(X_t, A_t) - J^*(\lambda, X_t)]\right\}\right]$$

$$\geq e^{\lambda\alpha[C(x,a) - J^*(\lambda,x)]}p_{xz}(a) + e^{\lambda\alpha[C(x,a) - J^*(\lambda,x)]}\sum_{y\neq z} p_{xy}(a)e^{\lambda h(y)},$$

where the inequality is due to the inclusion $\tilde{\delta} \in \mathcal{P}^*$; see (4.1). Recalling that $a \in B^*(x)$, this leads to

$$E_x^\pi\left[\exp\left\{\lambda\alpha \sum_{t=0}^{T-1}[C(X_t, A_t) - J^*(\lambda, X_t)]\right\}\Big| A_0 = a\right]$$

$$\geq \min_{b\in B^*(x)}\left[e^{\lambda\alpha[C(x,b) - J^*(\lambda,x)]}p_{xz}(b)\right.$$

$$\left. + e^{\lambda\alpha[C(x,b) - J^*(\lambda,x)]}\sum_{y\neq z} p_{xy}(b)e^{\lambda h(y)}\right],$$

and then, since this inequality holds for every action $a$ satisfying that $P_x^\pi[A_0 = a] > 0$,

$$E_x^\pi\left[\exp\left\{\lambda\alpha\sum_{t=0}^{T-1}[C(X_t, A_t) - J^*(\lambda, X_t)]\right\}\right]$$

$$\geq \min_{b\in B^*(x)}\left[e^{\lambda\alpha[C(x,b)-J^*(\lambda,x)]}p_{xz}(b)\right.$$

$$\left. + e^{\lambda\alpha[C(x,b)-J^*(\lambda,x)]}\sum_{y\neq z}p_{xy}(b)e^{\lambda h(y)}\right].$$

Using that $\pi \in \mathcal{P}^*$ and $x \in S$ are arbitrary, via (4.1), this inequality yields

(7.2)
$$e^{\lambda h(x)} \geq \min_{b\in B^*(x)}\left[e^{\lambda\alpha[C(x,b)-J^*(\lambda,x)]}p_{xz}(b)\right.$$

$$\left. + e^{\lambda\alpha[C(x,b)-J^*(\lambda,x)]}\sum_{y\neq z}p_{xy}(b)e^{\lambda h(y)}\right], \qquad x \in S.$$

On the other hand, by Theorem 6.1
$$e^{\lambda h(z)} \leq 1$$

which, combined with (7.2), implies that for every $x \in S$,

$$e^{\lambda h(x)} \geq \min_{b\in B^*(x)}\left[e^{\lambda\alpha[C(x,b)-J^*(\lambda,x)]}\sum_y p_{xy}(b)e^{\lambda h(y)}\right],$$

and then, multiplying both sides of this inequality by $e^{\lambda g(x)} = e^{\lambda[\alpha J^*(\lambda,x)+(1-\alpha)\|C\|]}$,

$$e^{\lambda g(x)+\lambda h(x)} \geq \min_{b\in B^*(x)}\left[e^{\lambda\alpha C(x,b)+(1-\alpha)\|C\|}\sum_y p_{xy}(b)e^{\lambda h(y)}\right];$$

since $\alpha C(x,b) + (1-\alpha)\|C\| \geq C(x,b)$, this yields that

$$e^{\lambda g(x)+\lambda h(x)} \geq \min_{b\in B^*(x)}\left[e^{\lambda C(x,b)}\sum_y p_{xy}(b)e^{\lambda h(y)}\right], \qquad x \in S.$$

Therefore, the pair $(g(\cdot), h(\cdot))$ satisfies the second condition of Definition 3.2, and it follows that

$$\alpha J^*(\lambda, \cdot) + (1-\alpha)\|C\| \in \mathcal{G}.$$

This inclusion is valid for each $\alpha \in (0,1)$, so that

$$J^*(\lambda, x) \geq \inf_{g\in\mathcal{G}} g(x), \qquad x \in S,$$

and, via Lemma 3.4, this implies that $J^*(\lambda, x) = \inf_{g\in\mathcal{G}} g(x)$ for every state $x$, completing the proof of Theorem 3.5. $\square$

REMARK 7.2. As a consequence of the results presented in this paper, two main problems remain open:

(i) Find (nontrivial) conditions under which the optimal value function $J(\lambda,\cdot)$ belongs to set $\mathcal{G}$ and there exists a solution to the dynamic programming equation.

(ii) Find an efficient algorithm to approximate the optimal value function and obtain $\varepsilon$-optimal stationary policies.

**Acknowledgment.** The authors are deeply grateful to the referee and the Associate Editor for their thoughtful reading of the original manuscript and helpful suggestions to improve the paper.

DEPARTAMENTO DE ESTADÍSTICA Y CÁLCULO
UNIVERSIDAD AUTÓNOMA AGRARIA
  ANTONIO NARRO
BUENAVISTA, SALTILLO COAH 25315
MEXICO
E-MAIL: rcavazos@narro.uaaan.mx

CENTRO DE INVESTIGACIÓN
  EN MATEMÁTICAS
APARTADO POSTAL 402
GUANAJUATO, GTO 36000
MEXICO
E-MAIL: dher@cimat.mx